\documentclass{imsart}
\usepackage[utf8]{inputenc}
\usepackage[T1]{fontenc}

\usepackage{amsmath,amsbsy,amsthm,amssymb}
\usepackage{multirow,multicol,tabularx,booktabs}
\usepackage{graphicx,placeins,url}
\usepackage{bbm,bm}
\usepackage[ruled]{algorithm2e}
\usepackage[shortlabels]{enumitem}
\usepackage[dvipsnames]{xcolor}
\usepackage[hidelinks]{hyperref}

\usepackage[top=2.5cm,bottom=2.5cm,right=2.5cm,left=2.5cm]{geometry}

\newcolumntype{L}{l<{\hspace{0.4cm}}}

\newtheorem{conj}{Conjecture}

\newtheorem{prop}[conj]{\bf Proposition}
\newtheorem{lemma}[conj]{\bf Lemma}
\newtheorem{rem}[conj]{\bf Remark}

\newtheorem{assumpt}{\bf Assumption}[section]

\def\bar{\overline}
\def\to{\rightarrow}

\def\indistrto{\buildrel {D} \over \longrightarrow}

\def\Ec{\mbox{$\mathcal E$}}
\def\Fc{\mbox{$\mathcal F$}}

\def\Hc{\mbox{$\mathcal H$}}

\def\Nc{\mbox{$\mathcal N$}}

\def\Sc{\mbox{$\mathcal S$}}

\def\Vc{\mbox{$\mathcal V$}}

\def\Xc{\mbox{$\mathcal X$}}

\def\Zc{\mbox{$\mathcal Z$}}

\def\EE{{\mathbb E}}

\def\HH{{\mathbb H}}

\def\NN{{\mathbb N}}

\def\Rb{\mbox{$\mathbb R$}}

\def\PP{ {\rm I} \kern-.15em {\rm P} }

\def\e{ {\bf e}}
\def\t{ {\bf t}}
\def\u{ {\bf u}}
\def\v{ {\bf v}}
\def\x{ {\bf x}}
\def\y{ {\bf y}}
\def\z{ {\bf z}}
\def\w{ {\bf w}}

\def\X{ {\bf X}}

\def\Z{ {\bf Z}}

\def\eps{\varepsilon}

\def\mds{\medskip}
\def \1{\mathbbm{1} }

\begin{document}

\begin{frontmatter}

\title{On kernel-based estimation of conditional Kendall's tau: finite-distance bounds and asymptotic behavior}
\runtitle{On kernel-based estimation of conditional Kendall's tau}

\begin{aug}
\author{\fnms{Alexis} \snm{Derumigny}\corref{}\ead[label=e1]{alexis.derumigny@ensae.fr}}
\and
\author{\fnms{Jean-David} \snm{Fermanian}\ead[label=e2]{jean-david.fermanian@ensae.fr}}

\runauthor{A. Derumigny and J.-D. Fermanian}

\affiliation{CREST-ENSAE}

\address{CREST-ENSAE, \\
5, avenue Henry Le Chatelier \\
91764 Palaiseau cedex, France. \\
\printead{e1}, 
\protect\printead*{e2}}
\end{aug}



\begin{abstract}
We study nonparametric estimators of conditional Kendall's tau, a measure of concordance between two random variables given some covariates. We prove non-asymptotic pointwise and uniform bounds, that hold with high probabilities. We provide ``direct proofs'' of the consistency and the asymptotic law of conditional Kendall's tau. A simulation study evaluates the numerical performance of such nonparametric estimators.    
\end{abstract}








\end{frontmatter}

\section{Introduction}

In the field of dependence modeling, it is common to work with dependence measures. 
Contrary to usual linear correlations, most of them have the advantage of being defined without any condition on moments, and of being invariant to changes in the underlying marginal distributions.
Such summaries of information are very popular and can be explicitly written as functionals of the underlying copulas: Kendall's tau, Spearman's rho, Blomqvist's coefficient...
See Nelsen~\cite{nelsen2007introduction} for an introduction.
In particular, for more  than a century (Spearman (1904), Kendall (1938)), Kendall's tau has become a popular dependence measure in $[-1,1]$.
It quantifies the positive or negative dependence between two random variables $X_1$ and $X_2$. 
Denoting by $C_{1,2}$ the unique underlying copula of $(X_1,X_2)$ that are assumed to be continuous, 
their Kendall's tau can be directly defined as
\begin{align}
    \tau_{1,2}
    &:= 4 \int_{[0,1]^2} C_{1,2}(u_1, u_2) \,
    C_{1,2}(du_1, du_2) - 1 \label{def_tau_cop}\\
    &= \PP \big( (X_{1,1}-X_{2,1})(X_{1,2}-X_{2,2}) > 0 \big)
    - \PP \big( (X_{1,1}-X_{2,1})(X_{1,2}-X_{2,2}) < 0 \big), \nonumber
\end{align}
where $(X_{i,1},X_{i,2})_{i=1,2}$ are two independent versions of $\X:=(X_1, X_2)$. This measure is then interpreted as the probability of observing a \textit{concordant pair} minus the probability of observing a \textit{discordant pair}.
See~\cite{kruskal1958} for an historical perspective on Kendall's tau.
Its inference is discussed in many textbooks (see~\cite{hollander1973} or~\cite{lehmann1975}, e.g.).
Its links with copulas and other dependence measures can be found in~\cite{nelsen2007introduction} or~\cite{joe1997multivariate}.

\mds

Similar dependence measure can be introduced in a conditional setup, when a $p$-dimensional covariate $\Z$ is available. 
When hundreds of papers refer to Kendall's tau, only a few of them have considered conditional Kendall's tau (as defined below) until now.
The goal is now to model the dependence between the two components $X_1$ and $X_2$, given the vector of covariates $\Z$. Logically, we can invoke the conditional copula $C_{1,2|\Z=\z}$ of $(X_1, X_2)$ given $\Z = \z$ for any point $\z \in \Rb^p$ (see Patton \cite{patton2006estimation, patton2006modelling}), and the corresponding conditional Kendall's tau would be simply defined as
\begin{align*}
    \tau_{1,2|\Z=\z}
    &:= 4 \int_{[0,1]^2} C_{1,2|\Z=\z}(u_1, u_2) \,
    C_{1,2|\Z=\z}(du_1, du_2) - 1 \\
    &= \PP \big( (X_{1,1}-X_{2,1})(X_{1,2}-X_{2,2}) > 0
    \big| \Z_1 = \Z_2 = \z \big) \\
    &\hspace{4cm} - \PP \big( (X_{1,1}-X_{2,1})(X_{1,2}-X_{2,2}) < 0
    \big| \Z_1 = \Z_2 = \z \big),
\end{align*}
where $(X_{i,1},X_{i,2},\Z_i)_{i=1,2}$ are two independent versions of $(X_1, X_2, \Z)$. As above, this is the probability of observing a \textit{concordant pair} minus the probability of observing a \textit{discordant pair}, conditionally on $\Z_1$ and $\Z_2$ being both equal to $\z$.
Note that, as conditional copulas themselves, conditional Kendall's taus are invariant w.r.t. increasing transformations of the conditional margins $X_1$ and $X_2$, given $\Z$. Of course, if $\Z$ is independent of $(X_1, X_2)$ then, for every $\z \in \Rb^p$, the conditional Kendall's tau $\tau_{1,2|\Z = \z}$ is equal to the (unconditional) Kendall's tau $\tau_{1,2}$.
\mds

Conditional Kendall's tau, and more generally conditional dependence measures, are of interest per se because they allow to summarize the evolution of the dependence between $X_1$ and $X_2$, when the covariate $\Z$ is changing. 
Surprisingly, their nonparametric estimates have been introduced in the literature only a few years ago (\cite{gijbels2011conditional},\cite{veraverbeke2011ScandinJ},\cite{fermanian2012time}) and their properties have not yet been fully studied in depth.
Indeed, until now and to the best of our knowledge, the theoretical properties of nonparametric conditional Kendall's tau estimates have been obtained ``in passing'' in the literature, as a sub-product of the weak-convergence of conditional copula processes (\cite{veraverbeke2011ScandinJ}) or as intermediate quantities that will be ``plugged-in'' (\cite{fermanian2015single}). Therefore, such properties have been stated under too demanding assumptions. In particular, some assumptions were related to the estimation of conditional margins, while this is not required because Kendall's tau are based on ranks. In this paper, we directly study nonparametric estimates $\hat\tau_{1,2|\z}$ without relying on the theory/inference of copulas. Therefore, we will state their main usual statistical properties: exponential bounds in probability, consistency, asymptotic normality. 

\mds

Our $\tau_{1,2|\Z=\z}$ has not to be confused with the so-called ``conditional Kendall's tau'' in the case of truncated data (\cite{tsai1990truncation},~\cite{martin2005testing}), in the case of semi-competing risk models (\cite{lakhal2008estimating},~\cite{hsieh2015nonparametric}), or for other partial information schemes (~\cite{chaieb2006estimating},~\cite{kim2015estimation}, among others). Indeed, particularly in biostatistics or reliability, the inference of dependence models under truncation/censoring can be led by considering some types of conditional Kendall's tau, given some algebraic relationships among the underlying random variables. This would induce conditioning by subsets. At the opposite, we will consider only pointwise conditioning events in this paper, under a nonparametric point-of-view. 
Nonetheless, such pointwise events can be found in the literature, but in some parametric or semi-parametric particular frameworks, as for the identifiability of frailty distributions in bivariate proportional models (~\cite{oakes1989bivariate},~\cite{manatunga1996measure}). 
Other related papers are~\cite{asimit2016tail} or~\cite{liu2015interval}, that are dealing with extreme co-movements (bivariate extreme-value theory). There, the tail conditioning events of Kendall's tau have probabilities that go to zero with the sample size.

\mds

In Section \ref{section:presentation_estimators}, different kernel-based estimators of the conditional Kendall's tau are discussed. In Section~\ref{section:theoretical_results}, the theoretical properties of the latter estimators are proved, first with finite-distance bounds and then under an asymptotic point-of-view.
A short simulation study is provided in Section~\ref{section:simulation_study}.
Proofs are postponed into the appendix.

\mds

\section{Definition of several kernel-based estimators of \texorpdfstring{$\tau_{1,2|\z}$}{tau12 cond z}}
\label{section:presentation_estimators}

Let $(X_{i,1}, X_{i,2}, \Z_i),$ $i=1,\dots, n$ be an i.i.d. sample distributed as $(X_1, X_2, \Z)$, and $n \geq 2$.
Assuming continuous underlying distributions, there are several equivalent ways of defining the conditional Kendall's tau:
\begin{align*}
    \tau_{1,2|\Z=\z}
    &= 4 \, \PP \big( X_{1,1} > X_{2,1}, X_{1,2} > X_{2,2}
    \big| \Z_1 = \Z_2 = \z \big) - 1  \\
    &= 1 - 4 \, \PP \big( X_{1,1} > X_{2,1}, X_{1,2} < X_{2,2}
    \big| \Z_1 = \Z_2 = \z \big)
    \displaybreak[0] \\
    &= \PP \big( (X_{1,1}-X_{2,1})(X_{1,2}-X_{2,2}) > 0
    \big| \Z_1 = \Z_2 = \z \big) 
    \nonumber \\ & \hspace{2cm}
    - \PP \big( (X_{1,1}-X_{2,1})(X_{1,2}-X_{2,2}) < 0
    \big| \Z_1 = \Z_2 = \z \big).
\end{align*}
Motivated by each of the latter expressions, we introduce several kernel-based estimators of $\tau_{1,2|\Z=\z}$: 
\begin{align*}
    \hat \tau_{1,2|\Z=\z}^{(1)}
    &:= 4 \sum_{i=1}^n \sum_{j=1}^n w_{i,n}(\z) w_{j,n}(\z)
    \1 \big\{ X_{i,1} < X_{j,1} , X_{i,2} < X_{j,2} \big\} - 1, \displaybreak[0] \\
    \hat \tau_{1,2|\Z=\z}^{(2)}
    &:= \sum_{i=1}^n \sum_{j=1}^n w_{i,n}(\z) w_{j,n}(\z)
    \Big( \1 \big\{ (X_{i,1} - X_{j,1}). (X_{i,2} - X_{j,2}) > 0 \big\} \\
    &\hspace{6cm} - \1 \big\{ (X_{i,1} - X_{j,1}). (X_{i,2} - X_{j,2}) < 0 \big\} \Big), \\
    \hat \tau_{1,2|\Z=\z}^{(3)}
    &:= 1 - 4 \sum_{i=1}^n \sum_{j=1}^n w_{i,n}(\z) w_{j,n}(\z)
    \1 \big\{ X_{i,1} < X_{j,1} , X_{i,2} > X_{j,2} \big\},
\end{align*}
where $\1$ denotes the indicator function, $w_{i,n}$ is a sequence of weights given by
\begin{equation}
    w_{i,n}(\z) = \frac{K_h(\Z_i-\z)}
    {\sum_{j=1}^n K_h(\Z_j-\z) },
    \label{def:weights_w_in}
\end{equation}
with $K_h(\cdot):= h^{-p} K(\cdot/h)$ for some kernel $K$ on $\Rb^p$, and $h=h(n)$ denotes a usual bandwidth sequence that tends to zero when $n\to\infty$.
In this paper, we have chosen usual Nadaraya-Watson weights. 
Obviously, there are alternatives (local linear, Priestley-Chao, Gasser-M\"uller, etc., weight), that would lead to different theoretical results. 

\mds

The estimators $\hat \tau_{1,2|\Z=\z}^{(1)},$
$\hat \tau_{1,2|\Z=\z}^{(2)}$ and
$\hat \tau_{1,2|\Z=\z}^{(3)}$ look similar, but they are nevertheless different, as shown in Proposition~\ref{prop:relationship_hat_tau_i}.
These differences are due to the fact that all the $\hat \tau_{1,2|\Z=\z}^{(k)},$ $k=1, 2, 3$ are affine transformations of a double-indexed sum, on every pair $(i,j)$, including the diagonal terms where $i=j$.
The treatment of these diagonal terms is different for each of the three estimators defined above.
Indeed, setting $s_n:= \sum_{i=1}^n w_{i,n}^2(\z),$
it can be easily proved that
$\hat \tau_{1,2|\Z=\z}^{(1)}$ takes values in the interval
$[- 1 \, , \, 1 - 2 s_n ]$,
$\hat \tau_{1,2|\Z=\z}^{(2)}$ in
$[- 1 + s_n \, , \, 1 - s_n ]$,
and $\hat \tau_{1,2|\Z=\z}^{(3)}$ in 
$[- 1 + 2 s_n \, , \, 1 ]$.
Moreover, there exists a direct relationship between these estimators, given by the following proposition.
\begin{prop}
    Almost surely, $\hat \tau_{1,2|\Z=\z}^{(1)} + s_n = \hat \tau_{1,2|\Z=\z}^{(2)} = \hat \tau_{1,2|\Z=\z}^{(3)} - s_n$, where $s_n:= \sum_{i=1}^n w_{i,n}^2(\z)$.
    \label{prop:relationship_hat_tau_i}
\end{prop}
This proposition is proved in~\ref{proof:prop:relationship_hat_tau_i}.
As a consequence, we can easily rescale the previous estimators so that the new estimator will take values in the whole interval $[-1, 1]$.
This would yield
\begin{align*}
    \tilde \tau_{1,2|\Z=\z}
    := \frac{\hat\tau_{1,2|\Z=\z}^{(1)}}{1-s_n} + \frac{s_n}{1-s_n}
    = \frac{\hat\tau_{1,2|\Z=\z}^{(2)}}{1-s_n}
    = \frac{\hat\tau_{1,2|\Z=\z}^{(3)}}{1-s_n} - \frac{s_n}{1-s_n}\cdot
\end{align*}

\mds

Note that none of the latter estimators depends on any estimation of conditional marginal distributions. In other words, we only have to conveniently choose the weights $w_{i,n}$ to obtain an estimator of the conditional Kendall's tau.
This is coherent with the fact that conditional Kendall's taus are invariant with respect to conditional marginal distributions.
Moreover, note that, in the definition of our estimators, the inequalities are strict (there are no terms corresponding
to the cases $i=j$). This is inline with the
definition of (conditional) Kendall's tau itself through concordant/discordant pairs of observations.


\mds

The definition of $\hat \tau_{1,2|\Z=\z}^{(1)}$ can be motivated as follows.
For $j=1,2$, let $\hat F_{j|\Z}(\cdot | \Z=\z)$ be an estimator of the conditional cdf of $X_j$ given $\Z=\z$.
Then, a usual estimator of the conditional copula of $X_1$ and $X_2$ given $\Z=\z$ is
\begin{equation*}
    \hat C_{1,2|\Z}(u_1,u_2 | \Z=\z) := \sum_{i=1}^n w_{i,n}(\z)
    \1 \big\{ \hat F_{1|\Z}(X_{i,1} | \Z=\z) \leq u_1 \, , \,
    \hat F_{2|\Z}(X_{i,2} | \Z=\z) \leq  u_2 \big\}.
\end{equation*}
See~\cite{veraverbeke2011ScandinJ} or~\cite{fermanian2012time}, e.g.
The latter estimator of the conditional copula can be plugged into~(\ref{def_tau_cop}) to define an estimator of the conditional Kendall's tau itself:
\begin{align}
    \hat \tau_{1,2|\Z=\z}
    &:= 4 \int \hat C_{1,2|\Z}(u_1, u_2 | \Z=\z) \,
    \hat C_{1,2|\Z}(du_1, du_2 | \Z=\z) - 1
    \label{tau_estim_veravEtAl} \\
    &= 4 \sum_{j=1}^n w_{j,n}(\z) \hat C_{1,2|\Z} \big(
    \hat F_{1|\Z} (X_{j,1} | \Z=\z) ,
    \hat F_{2|\Z} (X_{j,2} | \Z=\z) \big| \Z=\z \big)
    - 1 \nonumber.
\end{align}
Since the  functions $\hat F_{j|\Z}(\cdot | \Z=\z)$ are non-decreasing, this reduces to
\begin{eqnarray*}
\lefteqn{    \hat \tau_{1,2|\Z=\z}
    = 4 \sum_{i=1}^n \sum_{j=1}^n w_{i,n}(\z) w_{j,n}(\z)
    \1 \big\{ X_{i,1} \leq X_{j,1} , X_{i,2} \leq  X_{j,2} \big\} - 1  }\\
    &=& 4 \sum_{i=1}^n \sum_{j=1}^n w_{i,n}(\z) w_{j,n}(\z)
    \1 \big\{ X_{i,1} < X_{j,1} , X_{i,2} <  X_{j,2} \big\} - 1 + o_P(1)
    = \hat \tau_{1,2|\Z=\z}^{(1)}+ o_P(1).
\end{eqnarray*}
Veraverbeke et al.~\cite{veraverbeke2011ScandinJ}, Subsection 3.2, introduced their estimator of $\tau_{1,2|\z}$ by~(\ref{tau_estim_veravEtAl}). By the functional Delta-Method, they deduced its asymptotic normality as a sub-product of the weak convergence of the process 
$\sqrt{nh}\big(\hat C_{1,2|\Z}(\cdot, \cdot | z) - C_{1,2|\Z}(\cdot, \cdot | z)\big)$ when $\Z$ is univariate. 
In our case, we will obtain more and stronger theoretical properties of $\hat \tau_{1,2|\Z=\z}^{(1)}$ under weaker conditions by a more direct analysis based on ranks. In particular, we will not require any regularity condition on the conditional marginal distributions, contrary to~\cite{veraverbeke2011ScandinJ}. Indeed, in the latter paper, it is required that $F_{j|\Z}(\cdot | \Z=\z)$ has to be two times continuously differentiable (assumption $(\tilde{R}3)$) and its inverse has to be continuous (assumption $(R1)$). This is not satisfied for some simple univariate cdf as $F_j(t)=t\1(t\in [0,1])/2+ \1(t\in (1,2])/2+t\1(t\in (2,4])/4 + \1(t>4)$, for instance.
Note that we could justify $\hat \tau_{1,2|\Z=\z}^{(3)}$ in a similar way by considering conditional survival copulas.

\mds

Let us define $g_1, g_2, g_3$ by
\begin{align*}
    g_1(\X_i, \X_j)
    &:= 4 \1 \big\{ X_{i,1} < X_{j,1} , X_{i,2} < X_{j,2} \big\} - 1, \\
    g_2(\X_i, \X_j)
    &:= \1 \big\{ (X_{i,1} - X_{j,1})\times  (X_{i,2} - X_{j,2}) > 0 \big\}
    - \1 \big\{ (X_{i,1} - X_{j,1})\times  (X_{i,2} - X_{j,2}) < 0 \big\}, \\
    g_3(\X_i, \X_j)
    &:= 1 - 4 \1 \big\{ X_{i,1} < X_{j,1} , X_{i,2} > X_{j,2} \big\},
\end{align*}
where, for $i=1,\dots,n$, we set $\X_i := (X_{i,1},X_{i,2})$.
Clearly, $\hat \tau^{(k)}_{1,2|\z}$ is a smoothed estimator of $\EE[g_k(\X_1,\X_2) | \Z_1=\Z_2=\z]$, $k=1,2,3$.

\mds

Note that such dependence measures are of interest for the purpose of estimating (conditional or unconditional) copula models too.
Indeed, several popular parametric families of copulas have a simple one-to-one mapping between their parameter and the associated Kendall's tau (or Spearman's rho): Gaussian, Student with a fixed degree of freedom, Clayton, Gumbel and Frank copulas, etc.
Then, assume for instance that the conditional copula $C_{1,2|\Z=\z}$ belongs is a Gaussian copula with a parameter $\rho(\z)$. Then, by estimating its conditional Kendall's tau $\tau_{1,2|\Z=\z}$, we get an estimate of the corresponding parameter $\rho(\z)$, and finally of the conditional copula itself.
See~\cite{sabeti2014additive}, e.g.

\mds

The choice of the bandwidth $h$ could be done in a data-driven way, following the general conditional U-statistics framework detailed in Dony and Mason~\cite[Section~2]{dony2008uniform}.
Indeed, for any $k \in \{1,2,3\}$ and $\z \in \Zc$, denote by $\hat \tau_{-(i,j),\, 1,2|\Z=\z}^{(h,\, k)}$
the estimator $\hat \tau_{1,2|\Z=\z}^{(k)}$ that is made with the smoothing parameter $h$ and our dataset, when the $i$-th and $j$-th observations have been removed.
As a consequence, the random function $\hat \tau_{-(i,j),\, 1,2|\Z=\cdot}^{(h,\, k)}$ is independent of $\big( (\X_i, \Z_i) , (\X_j, \Z_j) \big)$.
As usual with kernel methods, it would be tempting to propose $h$ as the minimizer of the cross-validation criterion
\begin{align*}
    CV_{DM}(h) := \frac{2}{n(n-1)} \sum_{i,j=1}^n
    \Big( g_k(\X_i, \X_j) - \hat \tau_{-(i,j),\, 1,2| 
    \Z=(\Z_i + \Z_j)/2 }^{(h,\, k)} \Big)^2 K_h(\Z_i - \Z_j),
\end{align*}
for $k=1,2,3$ or for $\tilde\tau_{1,2|\Z=\cdot}$. 
The latter criterion would be a ``naively localized'' version of the usual cross-validation method.
Unfortunately, we observe that the function 
$h\mapsto CV_{DM}(h)$ is most often decreasing in the range of realistic bandwidth values. 
If we remove the weight $K_h(\Z_i - \Z_j)$, then there is no reason why $g_k(\X_i, \X_j)$ should be equal to 
$\hat \tau_{-(i,j),\, 1,2| \Z=(\Z_i + \Z_j)/2 }^{(k)}$ (on average), and we are not interested in the prediction of concordance/discordance pairs for which the $Z_i$ and $Z_j$ are far apart.
Therefore, a modification of this criteria is necessary. We propose to separate the choice of $h$ for the terms $g_k(\X_i, \X_j) - \hat \tau_{-(i,j),\, 1,2| \Z=(\Z_i + \Z_j)/2 }^{(h,\, k)}$ and the selection of the ``convenient pairs'' of observations $(i,j)$.
This leads to the new criterion
\begin{align}
    CV_{\tilde h}(h)
    &:= \frac{2}{n(n-1)} \sum_{i,j=1}^n
    \Big( g_k(\X_i, \X_j) - \hat \tau_{-(i,j),\, 1,2| 
    \Z=(\Z_i + \Z_j)/2 }^{(h,\, k)} \Big)^2
    \tilde K_{\tilde h}(\Z_i - \Z_j),
    \label{eq:def:new_CV_h}
\end{align}
with a potentially different kernel $\tilde K$ and a new fixed tuning parameter $\tilde h$.
Even if more complex procedures are possible, we suggest to simply choose
$\tilde K(\z) := \1 \{ |\z|_\infty  \leq 1\}$ and to calibrate $\tilde h$ so that only a fraction of the pairs $(i,j)$ has non-zero weights. In practice, set
$\tilde h$ as the empirical quantile of $\big( \{ |\Z_i - \Z_j|_\infty : 1 \leq i < j \neq n\}$ of order $2 N_{pairs}/(n(n-1))$, where $N_{pairs}$ is the number of pairs we want to keep.

\section{Theoretical results}
\label{section:theoretical_results}

\subsection{Finite distance bounds}

Hereafter, we will consider the behavior of conditional Kendall's tau estimates given $\Z=\z$ belongs to some fixed open subset $\Zc$ in $\Rb^p$.    
For the moment, let us state an instrumental result that is of interest per se. 
Let $\hat f_{\Z}(\z) := n^{-1} \sum_{j=1}^n K_h(\Z_j-\z)$ be the usual kernel estimator of the density $f_\Z$ of the conditioning variable $\Z$. Note that the estimators $\hat \tau_{1,2|\Z=\z}^{(k)},$ $k=1, \dots, 3$ are well-behaved only whenever $\hat f_{\Z}(\z) > 0$. Denote the joint density of $(\X,\Z)$ by $f_{\X,\Z}$. In our study, we need some usual conditions of regularity.

\mds



\begin{assumpt}
    The kernel $K$ is bounded, and set $\| K \|_{\infty} =: C_{K}$. It is symmetrical and satisfies $\int K = 1$, $\int |K| <\infty$.
    This kernel is of order $\alpha$ for some integer $\alpha > 1$: for all $j = 1, \dots, \alpha -1$ and every indices $i_1,\ldots,i_j$ in $\{1,\ldots,p\}$,
        $\int K(\u)  u_{i_1} \dots u_{i_j} \; d\u = 0$. Moreover, $\EE[K_h(\Z-\z)]>0$ for every $\z\in\Zc$ and $h>0$. 
        Set $\tilde K(\cdot):= K^2(\cdot)/\int K^2$ and $\| \tilde K \|_{\infty} =: C_{\tilde K}$.
    \label{assumpt:kernel_integral}
\end{assumpt}

\begin{assumpt}
    $f_\Z$ is $\alpha$-times continuously differentiable on $\Zc$ and there exists a constant $C_{K,\alpha}>0$ s.t.,
    for all $\z \in \Zc$,
    $$\int |K|(\u)
    \sum_{i_1, \dots, i_{ \alpha } = 1}^{p}
    |u_{i_1} \dots u_{i_{\alpha }}| \,
    \sup_{t\in [0,1]}\big| \frac{ \partial^{\alpha } f_{\Z}}{ \partial z_{i_1} \dots  \partial z_{i_{\alpha }}} (\z+th\u) \big| \, d\u \leq  C_{K,\alpha}.$$
    Moreover, $C_{\tilde K,2}$ denotes a similar constant replacing $K$ by $\tilde K$ and $\alpha$ by two.
    \label{assumpt:f_Z_Holder}
\end{assumpt}
\vspace{-0.3cm}

\begin{assumpt}
    There exist two positive constants $f_{\Z, min}$ and $f_{\Z, max}$ such that,
    for every $\z \in \Zc$, $f_{\Z, min} \leq f_{\Z}(\z) \leq f_{\Z, max}$.
    \label{assumpt:f_Z_max}
\end{assumpt}

\begin{prop}
    Under Assumptions \ref{assumpt:kernel_integral}-\ref{assumpt:f_Z_max} and if
    $ C_{K, \alpha} h^{\alpha}  / \alpha  ! \,
    < f_{\Z, min}$, for any $\z\in \Zc$, the estimator $\hat f_{\Z}(\z)$ is strictly positive with a probability larger than
    $$1 - 2 \exp \Big( - n h^p \big( f_{\Z, min} - C_{K, \alpha} h^{ \alpha}/\alpha ! \big)^2
    \, / \, \big( 2 f_{\Z, max} \int K^2 + (2/3) C_K ( f_{\Z, min} -
    C_{K, \alpha} h^{ \alpha}/\alpha !) \big) \Big).$$
    \label{cor:probaTau_Z_valid}
\end{prop}
\vspace{-0.3cm}
The latter proposition is proved in~\ref{proof:lemma:bound_f_hat_f}.
It guarantees that our estimators $\hat\tau^{(k)}_{1,2|\z}$, $k=1,\ldots,3$, are well-behaved with a probability close to one.
The next regularity assumption is necessary to explicitly control the bias of $\hat \tau_{1,2|\Z=\z}$.    
\begin{assumpt}
    For every $\x \in \Rb^2$, $\z \mapsto f_{\X, \Z}(\x, \z)$ is differentiable on $\Zc$ almost everywhere up to the order $ \alpha $. For every $0 \leq k \leq  \alpha $
    and every $1 \leq i_1, \dots, i_{ \alpha } \leq p$, let
    \begin{equation*}
        \Hc_{k,\vec{\iota}}(\u,\v,\x_1,\x_2,\z):= \sup_{t \in [0,1]} \bigg|
        \frac{\partial^{k} f_{\X, \Z}}{\partial z_{i_1} \dots  \partial z_{i_k}}
        \Big( \x_1, \z + t h\u \Big)
        \frac{\partial^{ \alpha -k} f_{\X, \Z}}
        {\partial z_{i_{k+1}} \dots  \partial z_{i_{ \alpha }}}
        \Big( \x_2, \z + t h\v \Big)\bigg|,
    \end{equation*}
    denoting $\vec{\iota}=(i_1,\ldots,i_\alpha)$.
    Assume that $\Hc_{k,\vec{\iota}}(\u,\v,\x_1,\x_2,\z)$ is integrable and there exists a finite constant $C_{\X\Z, \alpha} > 0$ such that, for every $\z \in \Zc$ and every $h<1$,
    \begin{align*}
        \int |K|(\u) |K|(\v)
        \sum_{k=0}^{ \alpha } \binom{ \alpha }{k}
        \sum_{i_1, \dots, i_{ \alpha } = 1}^{p}
        \Hc_{k, \vec{\iota}}(\u,\v,\x_1,\x_2,\z)
        |u_{i_1} \dots u_{i_k} v_{i_{k+1}} \dots v_{i_{ \alpha }}|
        \, d\u \, d\v\, d\x_1\, d\x_2
    \end{align*}
is less than $C_{\X\Z, \alpha}$.
    \label{assumpt:f_XZ_Holder}
\end{assumpt}

The next three propositions state pointwise and uniform exponential inequalities for the estimators $\hat \tau^{(k)}_{1,2|\Z=\z}$, when $k=1,2,3$. 
They are proved in~\ref{proof:prop:exponential_bound_KendallsTau}. We will denote $c_1 := c_3 := 4$ and $c_2 := 2$.
\begin{prop}[Exponential bound with explicit constants]
    Under Assumptions \ref{assumpt:kernel_integral}-\ref{assumpt:f_XZ_Holder},
    for every $t>0$ such that
    $ C_{K, \alpha} h^{\alpha} /  \alpha ! + t \leq f_{\Z, min}/2$ and every $t'>0$, if $C_{\tilde K,2} h^{2} < f_{\z}(\z)$, we have
    \begin{align*}
        &\PP \Bigg( |\hat \tau_{1,2|\Z=\z}^{(k)} - \tau_{1,2|\Z=\z} |
        >   \frac{c_k }{f^2_{\z}(\z)}\Big( \frac{C_{\X\Z, \alpha}   h^\alpha}{ \alpha  !}+ \frac{3f_{\z}(\z)\int K^2}{2nh^p}+ t' \Big)\times 
        \bigg( 1+ \frac{16 f^2_{\Z}(\z)} {f_{\Z, min}^3} 
        \Big( \frac{  C_{K, \alpha} h^{\alpha}} { \alpha  !} + t \Big) \bigg) \Bigg) \\
        &  \leq 2 \exp \Big( - \frac{n h^p t^2}{2 f_{\Z, max} \int K^2 + (2/3) C_K t} \Big)
        + 2 \exp \Big( - \frac{(n-1) h^{2p} t'{}^2 }{4 f_{\Z, max}^2 (\int K^2)^2 + (8/3) C_K^2 t'} \Big) \\
        &+ 2 \exp \bigg( - \frac{n h^p (f_{\z}(\z) - C_{\tilde K,2}h^2)^2 }{8 f_{\Z, max} \int \tilde K^2 + 4 C_{\tilde K} (f_{\z}(\z) - C_{\tilde K,2}h^2)/3} \bigg),
    \end{align*}
    for any $\z\in\Zc$ and every $k=1,2,3$.
    \label{prop:exponential_bound_KendallsTau}
\end{prop}

Alternatively, we can apply Theorem 1 in Major~\cite{major2006estimate} instead of the Bernstein-type inequality that has been used in the proof of Proposition~\ref{prop:exponential_bound_KendallsTau}. 
\begin{prop}[Alternative exponential bound without explicit constants]
    Under Assumptions \ref{assumpt:kernel_integral}-\ref{assumpt:f_XZ_Holder}, for every $t>0$ such that
    $ C_{K, \alpha} h^{\alpha} /  \alpha ! + t \leq f_{\Z, min}/2$ and every $t'>0$ s.t. $ t' \leq 2 h^p  (\int K^2)^3 f^3_{\Z,max} /  C_K^4$,  
    there exist some universal constants $C_2$ and $\alpha_2$ s.t. 
    \begin{align*}
        &\PP \Bigg( |\hat \tau_{1,2|\Z=\z}^{(k)} - \tau_{1,2|\Z=\z} |
        >         \frac{c_k }{f^2_{\z}(\z)}\Big( \frac{C_{\X\Z, \alpha}   h^\alpha}{ \alpha  !}+ \frac{3f_{\z}(\z)\int K^2}{2nh^p}+ t' \Big)\times 
    \bigg( 1+ \frac{16 f^2_{\Z}(\z)} {f_{\Z, min}^3} 
    \Big( \frac{  C_{K, \alpha} h^{\alpha}} { \alpha  !} + t \Big) \bigg)
        \Bigg) \\
        &  \leq 2 \exp \Big( - \frac{n h^p t^2}{2 f_{\Z, max} \int K^2 + (2/3) C_K t} \Big)
        + 2 \exp \bigg( - \frac{n h^p (f_{\z}(\z) - C_{\tilde K,2}h^2)^2 }{8 f_{\Z, max} \int \tilde K^2 + 4 C_{\tilde K} (f_{\z}(\z) - C_{\tilde K,2}h^2)/3} \bigg) \\
        &+ 2 \exp \Big( \frac{nh^p t^2}{32  \int K^2 (\int |K|)^2 f^3_{\Z,max} + 8 C_K\int |K|f_{\Z,max} t/3} \Big) 
        + C_2 \exp \bigg( - \frac{\alpha_2 n h^p t' }{ 8 f_{\Z, max} (\int K^2)} \bigg) ,
        \end{align*}
        for any $\z\in\Zc$ and every $k=1,2,3$, if $ C_{\tilde K, 2} h^{2}  < f_{\Z}(\z)$ and $6 h^p \big( \int |K|\big)^2 f_{\z,max} < \int K^2$.
\label{prop:exponential_bound_KendallsTau_major}
\end{prop}

\begin{rem}
    In Propositions~\ref{cor:probaTau_Z_valid},~\ref{prop:exponential_bound_KendallsTau} and~\ref{prop:exponential_bound_KendallsTau_major}, when the support of $K$ is included in $[-c,c]^p$ for some $c > 0$, $f_{\Z,max}$ can be replaced by a local bound
    $\sup_{\tilde \z \in \Vc(\z,\epsilon)} f_{\Z}(\tilde\z)$, denoting by $\Vc(\z,\epsilon)$ a closed ball
    of center $\z$ and any radius $\epsilon>0$, when $h\, c < \epsilon$. 
\end{rem}

As a corollary, the two latter result yield the weak consistency of
$\hat\tau_{1,2|\Z=\z}^{(k)}$ for every $\z\in\Zc$, when $nh^{2p}\to \infty$ (choose the constants $t$ and $t'\sim h^p$ sufficiently small, in Proposition~\ref{prop:exponential_bound_KendallsTau_major}, e.g.). 

\mds

It is possible to obtain uniform bounds, by slightly strengthening our assumptions.
Note that this next result will be true if $n$ is sufficiently large, when Proposition~\ref{prop:exponential_bound_KendallsTau_major} was true for every $n$.

\begin{assumpt}
    The kernel $K$ is Lipschitz on $(\Zc,\|\cdot\|_\infty)$, with a constant $\lambda_K$ and $\Zc$ is a subset of an hypercube in $\Rb^p$ whose volume is denoted by $\Vc$.  Moreover, $K$ and $K^2$ are regular in the sense of~\cite{gine2002rates} or~\cite{einmahl2005uniform}. 
    \label{assumpt:assum_Zc_vol}
\end{assumpt}

\begin{prop}[Uniform exponential bound]
    Under the assumptions~\ref{assumpt:kernel_integral}-\ref{assumpt:assum_Zc_vol}, there exist some constants $L_K$ and $\bar C_{K}$ (resp. $L_{\tilde K}$ and $\bar C_{\tilde K}$) that depend only on the VC characteristics of $K$ (resp. $\tilde K$), s.t., for every $\mu \in (0,1)$ such that $\mu f_{\z,min} < C_{\X\Z, \alpha}   h^\alpha / \alpha ! + b_K \int K^2 f_{\Z,max}/C_K$, if 
     $f_{\Z,max} < \tilde C_{\X\Z, 2} h^{2}/2 +  b_{\tilde K} \int \tilde K^2 f_{\Z,max}/C_{\tilde K} $,
    \begin{eqnarray*}
    \lefteqn{ \PP \Bigg( \sup_{\z\in\Zc}|\hat \tau_{1,2|\Z=\z}^{(k)} - \tau_{1,2|\Z=\z} |
        > \frac{c_k}{f^2_{\z,min}(1-\mu)^2} 
        \bigg( \frac{C_{\X\Z, \alpha}   h^\alpha}
        {  \alpha  !} + \frac{3f_{\z,max}\int K^2}{2nh^p} + t \bigg) \Bigg) }\\
        &  \leq & 
        L_K \exp \big( - C_{f,K} nh^p \big(\mu f_{\z,min} - \frac{C_{\X\Z, \alpha}   h^\alpha}{ \alpha  !}    \big)^2  \big) \\
        &+ &     C_2 D \exp \bigg( - \frac{\alpha_2 n t h^{p}}{  8 f_{\Z, max} (\int K^2)} \bigg) 
        +   L_{\tilde K} \exp \big( - C_{f, \tilde K} nh^p (f_{\z,max} - \tilde C_{\X\Z, 2}  h^{2})^2/4 \big) \\
        &+&
        2 \exp\big( -\frac{A_2 n h^p t^2 C_K^{-4}}{16^2 A_1^2  \int K^2 f_{\z,max}^3 (\int |K|)^2}  \big) + 2\exp(-\frac{A_2 n h^p t}{16 C_K^2 A_1} ) ,
        \end{eqnarray*}
        for $n$ sufficiently large, $k=1,2,3$, and for every $t>0$ s.t. $ t \leq 2 h^p  (\int K^2)^3 f^3_{\Z,max} / C_K^4$, 
$$ -16 A_1 C_K^2 A_{\bar g}\int K^2 f_{\z,max}^3 (\int |K|)^2 \ln(
h^p \int K^2 f_{\z,max}^3 (\int |K|)^2)  < n^{1/2}h^{p/2} t,\;\text{and}$$
$$ n h^p t \geq \big(\int K^2 \big)f_{\z,max}M_2 (p+\beta)^{3/2} \log\Big(\frac{4C_K^2 }{h^p f_{\z,max}\int K^2} \Big),\; \beta = \max\big(0, \frac{\log D}{\log n}\big),\; D:=\lceil \Vc\big( \frac{4 C_K \lambda_K}{h}\big)^p \rceil, $$
for some universal constants $C_2,\alpha_2,M_2,A_1,A_2$ and a constant $A_{\bar g}$ that depends on $K$ and $f_{\z,max}$.
\label{prop:exponential_bound_KendallsTau_uniform}
\end{prop}
We have denoted $ C_{f,K} :=  \log(1+b_K/(4L_K)) / (L_K b_K f_{\z,max} \int K^2)$, for any arbitrarily chosen constant $b_K\geq \bar C_{K}$. 
Similarly, $  C_{f,\tilde K} :=  \log(1+b_{\tilde K}/(4L_{\tilde K})) / (L_{\tilde K} b_{\tilde K} f_{\z,max} \int \tilde K^2)$, $b_{\tilde K} \geq \bar C_{\tilde K}$.

\mds

\subsection{Asymptotic behavior}

The previous exponential inequalities are not optimal to prove usual asymptotic results. Indeed, they directly or indirectly rely on upper bounds of estimates, as in Hoeffding or Bernstein-type inequalities. In the case of kernel estimates, this implies the necessary condition $nh^{2p}\rightarrow \infty$, at least. By a direct approach, it is possible to state the consistency of $\hat \tau_{1,2|\Z=\z}^{(k)}$, $k=1,2,3$, and then of $\tilde \tau_{1,2|\Z=\z}$, under the weaker condition $nh^p\rightarrow \infty$. 
\begin{prop}[Consistency]
    Under Assumption~\ref{assumpt:kernel_integral}, if
    $n h_{n}^p \to \infty$,
    $\lim K(\t) | \t |^p = 0$ when $|\t | \to \infty$,
    $f_\Z$ and $\z\mapsto \tau_{1,2|\Z=\z}$ are continuous on $\Zc$, then $\hat \tau_{1,2|\Z=\z}^{(k)}$ tends to $\tau_{1,2|\Z=\z}$ in probability, when $n\to\infty$ for any $k=1,2,3$.
    \label{prop:consistency_hatTau}
\end{prop}
This property is proved in~\ref{proof:prop:consistency_hatTau}.
Moreover, Proposition~\ref{prop:exponential_bound_KendallsTau_uniform} does not allow to state the strong uniform consistency of $\hat \tau_{1,2|\Z=\z}^{(k)}$ because the threshold $t$ has to be of order $h^p$ at most. Here again, a direct approach is possible, nonetheless. 

\begin{prop}[Uniform consistency]
    Under Assumption~\ref{assumpt:kernel_integral}, assume that
    $n h_{n}^{2p} / \log n \to \infty$,
    $\lim K(\t) | \t |^p = 0$ when $|\t | \to \infty$,
    $K$ is Lipschitz,
    $f_\Z$ and $\z\mapsto \tau_{1,2|\Z=\z}$ are continuous on a bounded set $\Zc$, and there exists a lower bound 
    $f_{\Z, \min}$ s.t. $f_ {\Z,\min} \leq f_{\Z}(\z)$ for any $\z\in\Zc$. Then $\sup_{\z \in \Zc} \big| \hat \tau_{1,2|\Z=\z}^{(k)} - \tau_{1,2|\Z=\z} \big| \to 0$ almost surely, when $n\to\infty$ for any $k=1,2,3$.
    \label{prop:unif_consistency_hatTau}
\end{prop}

This property is proved in~\ref{proof:prop:unif_consistency_hatTau}.
To derive the asymptotic law of this estimator, we will assume:
\begin{assumpt}
    (i) $n h_{n}^p \to \infty$ and $n h_{n}^{p+2\alpha} \to 0$;
    (ii) $K(\, \cdot \,)$ is compactly supported.
    \label{assumpt:asymptNorm}
\end{assumpt}

\begin{prop}[Joint asymptotic normality at different points]
    Let $\z'_1, \dots, \z'_{n'}$ be fixed points in a set $\Zc \subset \Rb^{p}$.
    Assume~\ref{assumpt:kernel_integral},~\ref{assumpt:f_XZ_Holder},~\ref{assumpt:asymptNorm}, that the $\z'_i$ are distinct and that $f_\Z$ and
    $\z\mapsto f_{\X,\Z}(\x,\z)$ are continuous on $\Zc$,
    for every~$\x$.
    Then, as $ n \to \infty$, 
    $$(n h_{n}^p)^{1/2} \left( \hat \tau_{1,2|\Z=\z'_i} 
    - \tau_{1,2|\Z=\z'_i} \right)_{i=1, \dots, n'}
    \indistrto \Nc (0,  \HH^{(k)}),\; k=1,2,3,$$ 
    where $\hat \tau_{1,2|\Z=\z}$ denotes any of the estimators $\hat \tau_{1,2|\Z=\z}^{(k)}$, $k=1,2,3$ or $\tilde \tau_{1,2|\Z=\z}$, and $\HH$ is the $n' \times n'$ diagonal real matrix defined by
    \begin{align*}
        [\HH^{(k)}]_{i, j}
        = \frac{ 4 \int K^2 \1_{ \{ i = j \} }
        }{ f_\Z(\z'_i)} \big\{
        \EE[g_{k}(\X_1,\X) g_{k}(\X_2,\X)
        | \Z = \Z_1 = \Z_2 = \z'_{i}]
         - \tau_{1,2|\Z=\z'_i}^2
        \big\},
    \end{align*}
    for every $1 \leq i, j \leq n'$, and $(\X, \Z)$, $(\X_1, \Z_1)$, $(\X_2, \Z_2)$ are independent versions.
    \label{prop:asymptNorm_hatTau}
\end{prop}
This proposition is proved in~\ref{proof:prop:asymptNorm_hatTau}.

\begin{rem}
The latter results will provide some simple tests of the constancy of the function $\z\mapsto \tau_{1,2|\z}$, and then of the constancy of the associated conditional copula itself. This would test the famous ``simplifying assumption'' (``$\Hc_0:C_{1,2|\Z=\z}$ does not depend on the choice of $\z$''), a key assumption for vine modeling in particular: see~\cite{acar2012beyond} or~\cite{HobaekAasFrigessi} for a discussion,~\cite{derumigny2017tests} for a review and a presentation of formal tests for this hypothesis.
\end{rem}

\section{Simulation study}
\label{section:simulation_study}

In this simulation study, we draw i.i.d. random samples $(X_{i,1}, X_{i,2}, Z_i),$ $i=1, \dots, n$, 
with univariate explanatory variables ($p=1$).
We consider two settings, that correspond to bounded and/or unbounded explanatory variables respectively:
\begin{enumerate}
    \item $\Zc = ]0,1[$ and the law of $Z$ is uniform on $]0,1[$. Conditionally on $Z=z$, $X_1 | Z = z$ and $X_2 | Z = z$ both follow a Gaussian distribution $\Nc(z,1)$. Their associated conditional copula is Gaussian and their conditional Kendall's tau is given by $\tau_{1,2|Z=z} = 2z-1$.
    
    \item $\Zc = \Rb$ and the law of $Z$ is $\Nc(0,1)$. Conditionally on $Z=z$, $X_1 | Z = z$ and $X_2 | Z = z$ both follow a Gaussian distribution $\Nc(\Phi(z),1)$, where $\Phi(\cdot)$ is the cdf of the $\Z$. Their associated conditional copula is Gaussian and their conditional Kendall's tau is given by $\tau_{1,2|Z=z} = 2\Phi(z)-1$.
\end{enumerate}

These simple frameworks allow us to compare the numerical properties of our different estimators in different parts of the space, in particular when $Z$ is close to zero or one, i.e. when the conditional Kendall's tau is close to $-1$ or to $1$.
We compute the different estimators $\hat \tau_{1,2|\Z=\z}^{(k)}$ for $k=1, 2, 3$, and the symmetrically rescaled version $\tilde\tau_{1,2|z}$. 
The bandwidth $h$ is chosen as proportional to the usual ``rule-of-thumb'' for kernel density estimation, i.e. $h = \alpha_h  \hat \sigma(Z)  n^{-1/5}$ with $\alpha_h \in \{ 0.5, 0.75, 1, 1.5, 2\}$ and $n \in \{ 100, 500, 1000, 2000 \}$.
For each setting, we consider three local measures of goodness-of-fit:  for a given $z$ and for any 
Kendall's tau estimate (say $\hat \tau_{1,2|Z=z}$), let
\begin{itemize}
    \item the (local) bias: $Bias(z) := \EE[\hat \tau_{1,2| Z=z}] - \tau_{1,2| Z=z}$,
    \item the (local) standard deviation: $Sd(z) := \EE \Big[ \big( \hat \tau_{1,2| Z=z} - \EE[\hat \tau_{1,2| Z=z}] \big)^2 \Big]^{1/2}$,
    \item the (local) mean square-error: $MSE(z) := \EE \Big[ \big( \hat \tau_{1,2| Z=z} - \tau_{1,2| Z=z} \big)^2 \Big]$.
\end{itemize}
We also consider their integrated version w.r.t the usual Lebesgue measure on the whole support of $z$, respectively denoted by $IBias$, $ISd$ and $IMSE$.
Some results concerning these integrated measures are given in Table~\ref{tab:result_simulation} (resp. Table~\ref{tab:result_simulation_znorm}) for Setting $1$ (resp. Setting $2$), and for different choices of $\alpha_h$ and $n$.
For the sake of effective calculations of these measures, all the theoretical previous expectations are replaced by their empirical counterparts based on $500$ simulations.

\mds

For every $n$, the best results seem to be obtained with $\alpha_h = 1.5$ and the fourth (rescaled) estimator, particularly in terms of bias.
This is not so surprising, because the estimators $\hat\tau^{(k)}$, $k=1,2,3$, do not have the right support at a finite distance. 
Note that this comparative advantage of $\tilde\tau$ in terms of bias decreases with $n$, as expected. In terms of integrated variance, all the considered estimators behave more or less similarly, particularly when $n\geq 500$.  

\mds

\begin{figure}[p]
    \centering
    \includegraphics[height = 7cm]{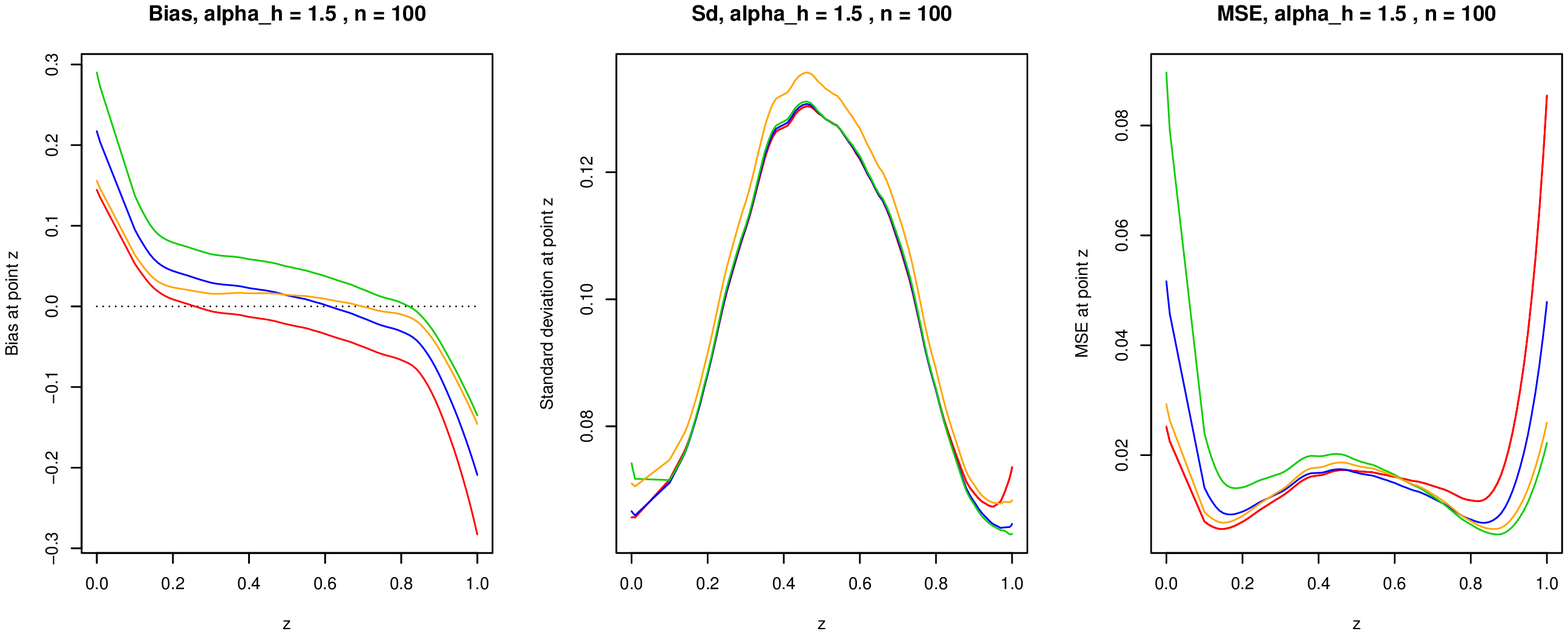}
    \caption{Local bias, standard deviation and MSE for the estimators $\hat \tau^{(1)}$ (red) , $\hat \tau^{(2)}$ (blue), $\hat \tau^{(3)}$~(green), $\tilde \tau$~(orange), with $n = 100$ and $\alpha_h = 1.5$ in Setting 1. The dotted line on the first figure is the reference at 0.}
    \label{fig:CKT_100}
\end{figure}

\begin{figure}
    \centering
    \includegraphics[height = 7cm]{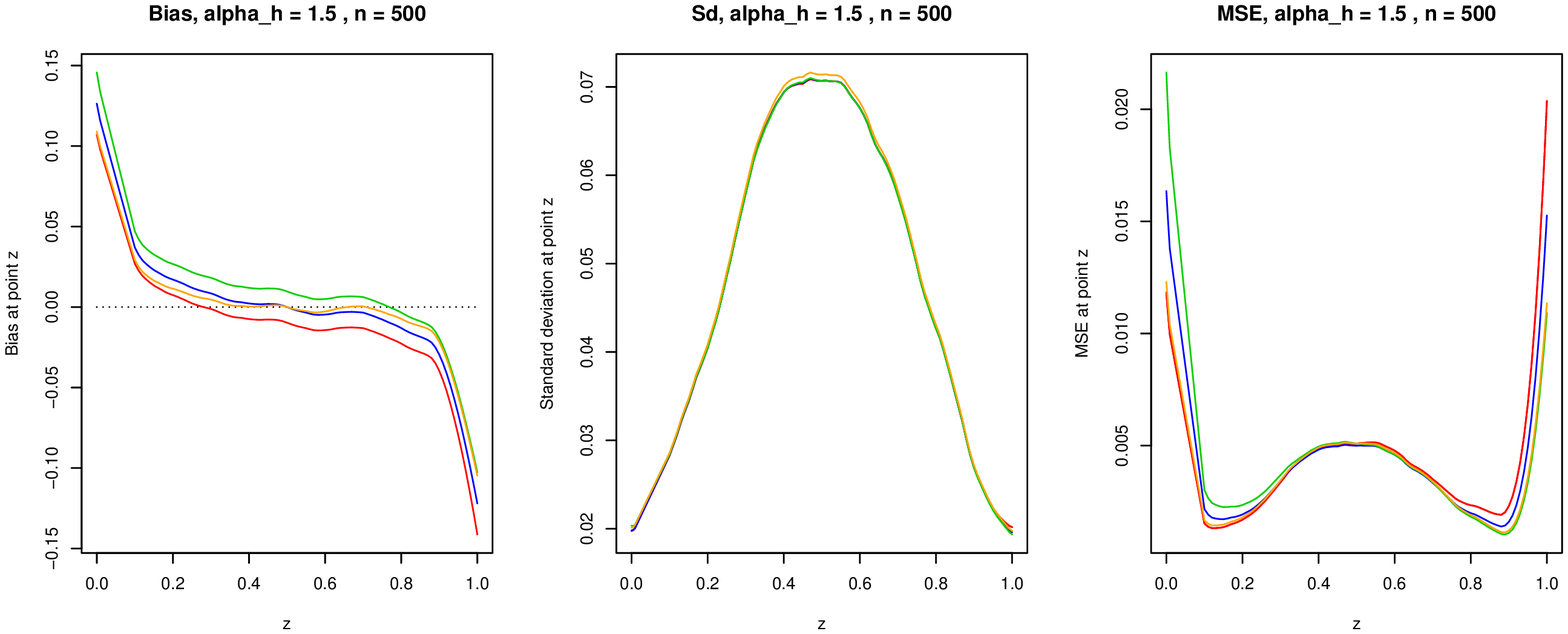}
    \caption{Local bias, standard deviation and MSE for the estimators $\hat \tau^{(1)}$ (red) , $\hat \tau^{(2)}$ (blue), $\hat \tau^{(3)}$~(green), $\tilde \tau$~(orange), with $n = 500$ and $\alpha_h = 1.5$ in Setting 1. The dotted line on the first figure is the reference at 0.}
    \label{fig:CKT_500}
\end{figure}

\begin{figure}
    \centering
    \includegraphics[height = 7cm]{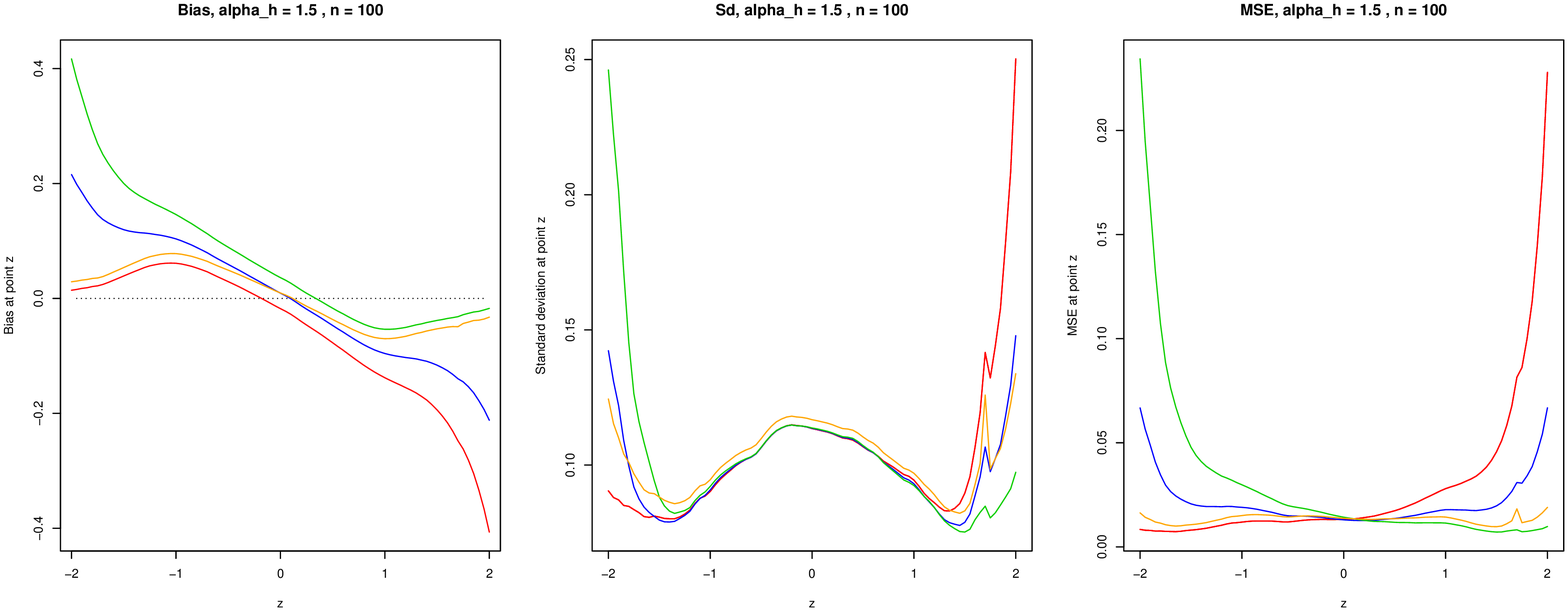}
    \caption{Local bias, standard deviation and MSE for the estimators $\hat \tau^{(1)}$ (red) , $\hat \tau^{(2)}$ (blue), $\hat \tau^{(3)}$~(green), $\tilde \tau$~(orange), with $n = 100$ and $\alpha_h = 1.5$ in Setting 2. The dotted line on the first figure is the reference at 0.}
    \label{fig:CKT_znorm_100}
\end{figure}

\begin{figure}
    \centering
    \includegraphics[height = 7cm]{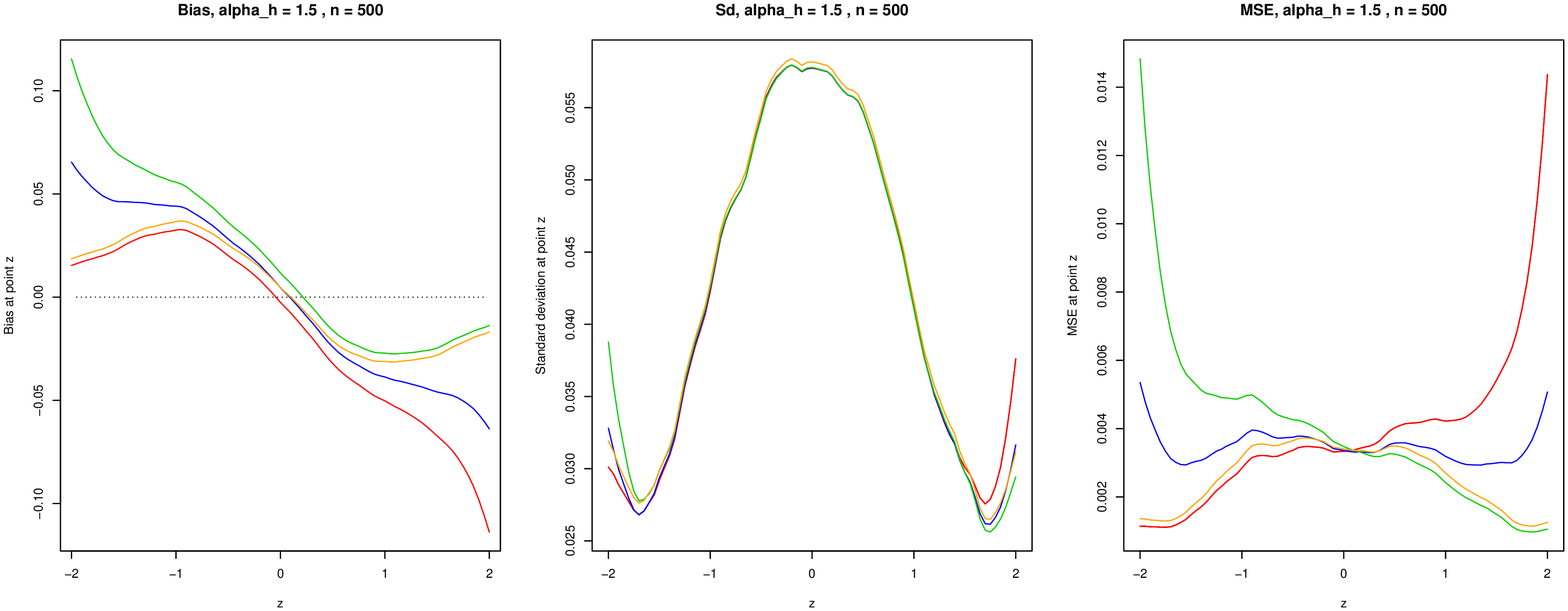}
    \caption{Local bias, standard deviation and MSE for the estimators $\hat \tau^{(1)}$ (red) , $\hat \tau^{(2)}$ (blue), $\hat \tau^{(3)}$~(green), $\tilde \tau$~(orange), with $n = 500$ and $\alpha_h = 1.5$ in Setting 2. The dotted line on the first figure is the reference at 0.}
    \label{fig:CKT_znorm_500}
\end{figure}

To illustrate our results for Setting 1 (resp. Setting 2), the functions $z \mapsto Bias(z)$, $Sd(z)$ and $MSE(z)$ have been plotted on Figures \ref{fig:CKT_100}-\ref{fig:CKT_500} (resp. Figures~\ref{fig:CKT_znorm_100}-\ref{fig:CKT_znorm_500}), both with our empirically optimal choice $\alpha_h = 1.5$.
We can note that, considering the bias, the estimator $\tilde \tau$ behaves similarly as $\hat \tau^{(1)}$ when the true $\tau$ is close to $-1$, and similarly as $\hat \tau^{(3)}$ when the true Kendall's tau is close to $1$. But globally, the best pointwise estimator is clearly obtained with the rescaled version $\tilde \tau_{1,2|\Z=\cdot}$, after a quick inspection of MSE levels, and even if the differences between our four estimators weaken for large sample sizes. The comparative advantage of $\tilde \tau_{1,2|\z}$ more clearly appears with Setting 2 than with Setting 1.
Indeed, in the former case, the support of $\Z$'s distribution is the whole line. Then $\hat f_{\Z}$ does not suffer any more from the boundary bias phenomenon, contrary to what happened with Setting 1. As a consequence, the biases induced by the definitions of $\hat \tau^{(k)}_{1,2|\z}$, $k=1,3$, appear more strinkingly in Figure~\ref{fig:CKT_znorm_100}, for instance: when $z$ is close to $(-1)$ (resp. $1$), the biases of $\hat \tau^{(1)}_{1,2|z}$ (resp. $\hat \tau^{(3)}_{1,2|z}$) and $\tilde \tau_{1,2|z}$ are close, when the bias $\hat \tau^{(3)}_{1,2|z}$ (resp. $\hat \tau^{(1)}_{1,2|z}$) is a lot larger.  
Since the squared biases are here significantly larger than the variances in the tails, $\tilde \tau_{1,2|z}$ provides the best estimator globally considering ''both sides'' together. But even in the center of $\Z$'s distribution, the latter estimator behaves very well. 

\mds

In Setting 2 where there is no boundary problem, we also try to estimate the conditional Kendall's tau using our cross-validation criterion~(\ref{eq:def:new_CV_h}), with $N_{pairs} = 1000$.
More precisely, denoting by $h^{CV}$ the minimizer of the cross-validation criterion, we try different choices $h=\alpha_h \times h^{CV}$ with $\alpha_h \in \{ 0.5, 0.75, 1, 1.5, 2\}$. The results in terms of integrated bias, standard deviation and MSE are given in Table~\ref{tab:result_simulation_znorm_hCV}. We do not find any substantial improvements compared to the previous Table~\ref{tab:result_simulation_znorm}, where the bandwidth was chosen ``roughly''. In Table~\ref{tab:summary_statistics_hCV}, we compare the average $h^{CV}$ with the previous choice of $h$. The expectation of $h^{CV}$ is always higher than the ``rule-of-thumb'' $h^{ref}$, but the difference between both decreases when the sample size $n$ increases. The standard deviation of $h^{CV}$ is quite high for low values of $n$, but decreases as a function of $n$. This may be seen as quite surprising given the fact that the number of pairs $N_{pairs}$ used in the computation of the criterion stays constant. Nevertheless, when the sample size increases, the selected pairs are better in the sense that the differences $|Z_i - Z_j|$ can become smaller as more replications of $Z_i$ are available.

\FloatBarrier

\begin{table}[p]
    \centering
    \renewcommand{\arraystretch}{1.3}
    \resizebox{\textwidth}{!}{%
    \begin{tabular}{cc ccc | ccc | ccc | ccc}
        &&
        \multicolumn{3}{c}{$n = 100$} &
        \multicolumn{3}{c}{$n = 500$} &
        \multicolumn{3}{c}{$n = 1000$} &
        \multicolumn{3}{c}{$n = 2000$} \\
        \cmidrule(lr){3-5} \cmidrule(lr){6-8}
        \cmidrule(lr){9-11} \cmidrule(lr){12-14}
        && IBias & ISd & IMSE & IBias & ISd & IMSE
        & IBias & ISd & IMSE & IBias & ISd & IMSE \\
        \parbox[t]{2mm}{\multirow{4}{*}{\rotatebox[origin=c]{90}{ 
        $\alpha_h = 0.5$}}} &  $\hat \tau_{1,2|\Z=\z}^{(1)}$  & -133 &  197 & 66.5 & -34.5 & 84.9 & 9.86 & -18.2 & 61.6 & 4.85 & -10.9 &  46 & 2.65 \\ 
         &  $\hat \tau_{1,2|\Z=\z}^{(2)}$  & -12.9 &  187 & 43.7 & -4.08 & 84.4 & 8.58 & -0.9 & 61.5 & 4.49 & -1.07 &   46 & 2.53 \\ 
         &  $\hat \tau_{1,2|\Z=\z}^{(3)}$  &  107 &  190 & 56.6 & 26.4 & 84.5 & 9.26 & 16.4 & 61.5 & 4.76 &  8.8 &   46 &  2.6 \\ 
         &  $\tilde \tau_{1,2|\Z=\z}$  & \textbf{-0.91} &  213 & 48.2 & \textbf{-1.18} & 86.9 & 8.55 & \textbf{0.733} & 62.4 & 4.46 & \textbf{-0.149} & 46.4 &  2.5 \\ 
        \hline\parbox[t]{2mm}{\multirow{4}{*}{\rotatebox[origin=c]{90}{ 
        $\alpha_h = 0.75$}}} &  $\hat \tau_{1,2|\Z=\z}^{(1)}$  &  -88 &  150 & 35.8 & -26.3 &   68 & 6.32 & -13.9 & 50.7 & 3.33 & -7.98 & 37.6 &  1.8 \\ 
         &  $\hat \tau_{1,2|\Z=\z}^{(2)}$  & -10.4 &  145 & 26.3 & -5.97 & 67.9 &  5.6 & -2.33 & 50.6 & 3.12 & -1.39 & 37.5 & 1.74 \\ 
         &  $\hat \tau_{1,2|\Z=\z}^{(3)}$  & 67.2 &  146 & 30.6 & 14.3 & 67.9 & 5.75 &  9.2 & 50.6 & 3.19 &  5.2 & 37.5 & 1.76 \\ 
         &  $\tilde \tau_{1,2|\Z=\z}$  & -2.06 &  157 & 26.7 & -3.99 & 69.2 & 5.49 & -1.21 & 51.2 & 3.05 & -0.76 & 37.8 & 1.69 \\
        \hline\parbox[t]{2mm}{\multirow{4}{*}{\rotatebox[origin=c]{90}{ 
        $\alpha_h = 1$}}} &  $\hat \tau_{1,2|\Z=\z}^{(1)}$  & -67.8 &  123 & 24.5 & -19.2 & 58.7 &  4.8 &  -11 & 43.1 & 2.52 & -6.34 &   33 & 1.44 \\ 
         &  $\hat \tau_{1,2|\Z=\z}^{(2)}$  & -9.99 &  121 &   19 & -3.95 & 58.6 & 4.39 & -2.35 & 43.1 & 2.39 & -1.39 &   33 &  1.4 \\ 
         &  $\hat \tau_{1,2|\Z=\z}^{(3)}$  & 47.8 &  122 & 20.9 & 11.3 & 58.7 & 4.47 & 6.34 & 43.1 & 2.41 & 3.57 &   33 & 1.41 \\ 
         &  $\tilde \tau_{1,2|\Z=\z}$  & -3.48 &  128 & 18.1 & -2.34 & 59.5 & 4.18 & -1.46 & 43.4 & 2.29 & -0.897 & 33.2 & 1.35 \\ 
        \hline\parbox[t]{2mm}{\multirow{4}{*}{\rotatebox[origin=c]{90}{ 
        $\alpha_h = 1.5$}}} &  $\hat \tau_{1,2|\Z=\z}^{(1)}$  & -44.6 &  101 & 17.5 & -15.9 & 50.4 & 4.12 & -9.7 & 35.9 & 2.13 & -5.52 & 27.6 & 1.28 \\ 
         &  $\hat \tau_{1,2|\Z=\z}^{(2)}$  & -5.81 &  100 & 14.9 & -5.68 & 50.3 & 3.84 & -3.84 & 35.9 & 2.02 & -2.18 & 27.6 & 1.24 \\ 
         &  $\hat \tau_{1,2|\Z=\z}^{(3)}$  &   33 &  101 & 15.5 & 4.58 & 50.3 & 3.77 & 2.01 & 35.9 & 1.99 & 1.15 & 27.6 & 1.23 \\ 
         &  $\tilde \tau_{1,2|\Z=\z}$  & -1.09 &  104 & \textbf{13.4} & -4.55 & 50.8 & \textbf{3.57} & -3.19 & 36.1 &  \textbf{1.9} & -1.83 & 27.7 & \textbf{1.18} \\ 
        \hline\parbox[t]{2mm}{\multirow{4}{*}{\rotatebox[origin=c]{90}{ 
        $\alpha_h = 2$}}} &  $\hat \tau_{1,2|\Z=\z}^{(1)}$  & -37.8 & \textrm{91.4} & 17.3 & -11.8 & \textrm{43.8} & 4.14 & -7.2 & \textrm{31.2} & 2.35 & -5.97 & \textrm{23.7} & 1.43 \\ 
        &  $\hat \tau_{1,2|\Z=\z}^{(2)}$  & -8.03 & \textrm{91.4} & 15.4 & -3.93 & \textrm{43.8} & 3.94 & -2.75 & \textrm{31.2} & 2.28 & -3.44 & \textrm{23.7} & 1.39 \\ 
        &  $\hat \tau_{1,2|\Z=\z}^{(3)}$  & 21.7 & 91.7 & 15.4 & 3.91 & \textrm{43.8} & 3.87 &  1.7 & \textrm{31.2} & 2.24 & -0.912 & \textrm{23.7} & 1.37 \\ 
        &  $\tilde \tau_{1,2|\Z=\z}$  & -4.5 & 94.2 & 13.5 & -3.01 & 44.1 & 3.62 & -2.24 & 31.3 & 2.12 & -3.16 & 23.8 & 1.32 \\ 
        \hline
        
    \end{tabular}
    }
    \caption{Results of the simulation in Setting 1. All values have been multiplied by 1000. Bold values indicate optimal choices for the chosen measure of performance.}
    \label{tab:result_simulation}
\end{table}

\begin{table}[p]
    \centering
    \renewcommand{\arraystretch}{1.3}
    \resizebox{\textwidth}{!}{%
    \begin{tabular}{cc ccc | ccc | ccc | ccc}
        &&
        \multicolumn{3}{c}{$n = 100$} &
        \multicolumn{3}{c}{$n = 500$} &
        \multicolumn{3}{c}{$n = 1000$} &
        \multicolumn{3}{c}{$n = 2000$} \\
        \cmidrule(lr){3-5} \cmidrule(lr){6-8}
        \cmidrule(lr){9-11} \cmidrule(lr){12-14}
        && IBias & ISd & IMSE & IBias & ISd & IMSE
        & IBias & ISd & IMSE & IBias & ISd & IMSE \\
        
        \parbox[t]{2mm}{\multirow{4}{*}{\rotatebox[origin=c]{90}{ 
        $\alpha_h = 0.5$}}}&  $\hat \tau_{1,2|\Z=\z}^{(1)}$  & -207 &  227 &  180 & -54.1 & 83.9 & 16.9 & -29.6 & 55.3 & 5.81 & -16.9 & 38.9 & 2.49 \\ 
        &  $\hat \tau_{1,2|\Z=\z}^{(2)}$  & \textbf{1.15} &  207 &   97 & \textbf{0.845} & 80.5 & 10.8 & 0.557 & 54.4 & 4.35 & \textbf{0.145} & 38.6 & 2.04 \\ 
        &  $\hat \tau_{1,2|\Z=\z}^{(3)}$  &  210 &  228 &  181 & 55.7 & 83.2 & 16.4 & 30.7 & 55.4 &  5.9 & 17.2 & 38.9 &  2.5 \\ 
        &  $\hat \tau_{1,2|\Z=\z}^{(4)}$  &  1.4 &  225 & 51.9 & 0.987 & 81.4 & 6.86 & \textbf{0.456} &   55 & 3.22 & 0.175 & 38.9 & 1.66 \\ 
        \hline\parbox[t]{2mm}{\multirow{4}{*}{\rotatebox[origin=c]{90}{ 
        $\alpha_h = 0.75$}}}&  $\hat \tau_{1,2|\Z=\z}^{(1)}$  & -144 &  175 & 98.6 & -33.3 & 60.6 &  7.5 & -19.8 & 41.9 & 3.12 & -10.6 & 30.5 & 1.42 \\ 
        &  $\hat \tau_{1,2|\Z=\z}^{(2)}$  & -2.33 &  163 & 56.2 & 1.73 & 59.4 & 5.56 & -0.0619 & 41.7 & 2.51 & 0.665 & 30.4 & 1.24 \\ 
        &  $\hat \tau_{1,2|\Z=\z}^{(3)}$  &  140 &  176 & 99.2 & 36.8 & 60.7 & 7.73 & 19.7 & 42.1 & 3.12 & 11.9 & 30.5 & 1.45 \\ 
        &  $\hat \tau_{1,2|\Z=\z}^{(4)}$  & -3.15 &  170 & 30.3 & 1.69 & 60.2 & 3.85 & -0.093 & 42.1 & 1.95 & 0.645 & 30.5 & 1.05 \\ 
        \hline\parbox[t]{2mm}{\multirow{4}{*}{\rotatebox[origin=c]{90}{ 
        $\alpha_h = 1$}}}&  $\hat \tau_{1,2|\Z=\z}^{(1)}$  & -99.8 &  143 & 57.7 & -24.9 & 50.9 & 5.06 & -13.5 & 36.6 & 2.28 & -6.92 & 26.6 & 1.09 \\ 
        &  $\hat \tau_{1,2|\Z=\z}^{(2)}$  & 1.17 &  132 & 34.6 & 0.903 & 50.4 & 4.02 & 1.16 & 36.5 & 1.97 & 1.46 & 26.6 & 0.994 \\ 
        &  $\hat \tau_{1,2|\Z=\z}^{(3)}$  &  102 &  139 & 54.4 & 26.7 &   51 & 5.13 & 15.8 & 36.6 & 2.33 & 9.83 & 26.6 & 1.11 \\ 
        &  $\hat \tau_{1,2|\Z=\z}^{(4)}$  & 2.51 &  138 & 20.1 & 0.897 & 50.9 & 2.89 & 1.16 & 36.7 & 1.56 & 1.48 & 26.7 & 0.847 \\ 
        \hline\parbox[t]{2mm}{\multirow{4}{*}{\rotatebox[origin=c]{90}{ 
        $\alpha_h = 1.5$}}}&  $\hat \tau_{1,2|\Z=\z}^{(1)}$  & -59.1 &  104 & 28.1 & -14.7 & 42.3 & 3.87 & -7.56 & 29.7 & 1.86 & -4.17 & 21.8 & 0.932 \\ 
        &  $\hat \tau_{1,2|\Z=\z}^{(2)}$  & 4.34 & 99.7 & 21.4 & 2.05 & 42.1 & 3.48 & 2.07 & 29.6 & 1.75 & 1.35 & 21.8 & 0.899 \\ 
        &  $\hat \tau_{1,2|\Z=\z}^{(3)}$  & 67.8 &  103 & 29.6 & 18.8 & 42.3 & 3.96 & 11.7 & 29.6 & 1.92 & 6.87 & 21.8 & 0.957 \\ 
        &  $\hat \tau_{1,2|\Z=\z}^{(4)}$  & 3.34 &  103 & \textbf{13.4} & 2.08 & 42.5 &  \textbf{2.6} & 2.08 & 29.7 & \textbf{1.39} & 1.35 & 21.8 & \textbf{0.755} \\ 
        \hline\parbox[t]{2mm}{\multirow{4}{*}{\rotatebox[origin=c]{90}{ 
        $\alpha_h = 2$}}}&  $\hat \tau_{1,2|\Z=\z}^{(1)}$  & -37.2 & 88.2 & 23.9 & -9.57 & 38.2 &  4.6 & -3.75 & 26.2 & 2.34 & -1.09 & \textrm{19.8} & 1.32 \\ 
        &  $\hat \tau_{1,2|\Z=\z}^{(2)}$  & 8.17 & \textrm{85.9} & 21.2 & 2.69 &   \textrm{38} & 4.45 & 3.32 & \textrm{26.1} &  2.3 & 2.99 & \textrm{19.8} & 1.32 \\ 
        &  $\hat \tau_{1,2|\Z=\z}^{(3)}$  & 53.5 & 87.4 & 25.3 & 14.9 & 38.1 & 4.74 & 10.4 & 26.2 & 2.41 & 7.08 & \textrm{19.8} & 1.36 \\ 
        &  $\hat \tau_{1,2|\Z=\z}^{(4)}$  & 8.47 & 88.5 &   15 & 2.69 & 38.4 & 3.59 & 3.33 & 26.3 & 1.93 &    3 & 19.9 & 1.15 \\ 
        \hline
        
    \end{tabular}
    }
    \caption{Results of the simulation in Setting 2. All values have been multiplied by 1000. Bold values indicate optimal choices for the chosen measure of performance.}
    \label{tab:result_simulation_znorm}
\end{table}

\begin{table}[p]
    \centering
    \renewcommand{\arraystretch}{1.3}
    \resizebox{\textwidth}{!}{%
    \begin{tabular}{cc ccc | ccc | ccc | ccc}
        &&
        \multicolumn{3}{c}{$n = 100$} &
        \multicolumn{3}{c}{$n = 500$} &
        \multicolumn{3}{c}{$n = 1000$} &
        \multicolumn{3}{c}{$n = 2000$} \\
        \cmidrule(lr){3-5} \cmidrule(lr){6-8}
        \cmidrule(lr){9-11} \cmidrule(lr){12-14}
        && IBias & ISd & IMSE & IBias & ISd & IMSE
        & IBias & ISd & IMSE & IBias & ISd & IMSE \\
        
    \parbox[t]{2mm}{\multirow{4}{*}{\rotatebox[origin=c]{90}{ 
    $\alpha_h = 0.5$}}}&  $\hat \tau_{1,2|\Z=\z}^{(1)}$  & -111 &  154 & 66.2 & -36.9 & 66.8 & 9.01 & -22.4 & 48.2 & 4.06 & -12.9 & 36.1 & 2.04 \\ 
    &  $\hat \tau_{1,2|\Z=\z}^{(2)}$  & \textbf{0.0488} &  137 & 36.3 & \textbf{0.236} & 64.2 & 6.45 & \textbf{0.546} & 46.8 & 3.14 & \textbf{1.29} & 35.7 & 1.78 \\ 
    &  $\hat \tau_{1,2|\Z=\z}^{(3)}$  &  111 &  151 & 60.6 & 37.4 & 66.3 & 8.88 & 23.5 & 47.2 & 4.07 & 15.5 & 36.2 & 2.18 \\ 
    &  $\hat \tau_{1,2|\Z=\z}^{(4)}$  & 1.38 &  132 & 18.3 & 0.27 & 64.5 & 4.49 & 0.61 & 46.8 & 2.36 & 1.29 & 35.6 & 1.49 \\ 
    \hline\parbox[t]{2mm}{\multirow{4}{*}{\rotatebox[origin=c]{90}{ 
    $\alpha_h = 0.75$}}}&  $\hat \tau_{1,2|\Z=\z}^{(1)}$  & -67.4 &  117 & 35.7 & -23.3 & 52.1 & 5.27 & -13.9 & 37.8 &  2.4 & -7.6 &   29 &  1.3 \\ 
    &  $\hat \tau_{1,2|\Z=\z}^{(2)}$  & 4.32 &  108 & 23.5 & 0.809 & 50.7 & 4.21 & 1.03 & 37.2 & 2.07 & 1.78 & 28.8 & 1.21 \\ 
    &  $\hat \tau_{1,2|\Z=\z}^{(3)}$  & 76.1 &  119 & 35.4 & 24.9 & 51.6 & 5.12 &   16 & 37.6 & 2.49 & 11.2 & 29.1 & 1.39 \\ 
    &  $\hat \tau_{1,2|\Z=\z}^{(4)}$  & 4.98 &  106 & \textbf{13.3} & 0.86 & 51.6 & 3.13 & 1.03 & 37.5 & 1.63 & 1.81 & 28.9 & 1.02 \\ 
    \hline\parbox[t]{2mm}{\multirow{4}{*}{\rotatebox[origin=c]{90}{ 
    $\alpha_h = 1$}}}&  $\hat \tau_{1,2|\Z=\z}^{(1)}$  &  -43 &  101 &   28 & -15.8 & 45.7 & 4.44 & -9.51 & 33.1 & 2.04 & -4.68 & 25.1 & 1.07 \\ 
    &  $\hat \tau_{1,2|\Z=\z}^{(2)}$  & 7.87 & 93.1 & 22.4 & 2.01 & 44.8 & 3.91 & 1.57 & 32.7 & 1.87 & 2.29 & 24.9 & 1.03 \\ 
    &  $\hat \tau_{1,2|\Z=\z}^{(3)}$  & 58.8 & 97.6 & 27.2 & 19.8 & 45.3 & 4.41 & 12.7 & 32.9 &  2.1 & 9.27 & 25.1 & 1.14 \\ 
    &  $\hat \tau_{1,2|\Z=\z}^{(4)}$  & 8.51 &   98 & 15.7 & 2.05 &   46 & \textbf{3.01} & 1.57 & 33.1 &  \textbf{1.5} & 2.33 & 25.1 & \textbf{0.871} \\ 
    \hline\parbox[t]{2mm}{\multirow{4}{*}{\rotatebox[origin=c]{90}{ 
    $\alpha_h = 1.5$}}}&  $\hat \tau_{1,2|\Z=\z}^{(1)}$  & -16.1 & 95.6 & 41.7 & -6.36 &   43 & 6.35 & -4.04 & 30.6 & 2.87 & -1.11 & 22.1 & 1.34 \\ 
    &  $\hat \tau_{1,2|\Z=\z}^{(2)}$  & 14.9 & \textbf{92.6} & 40.4 & 5.08 & 42.6 &  6.2 & 3.17 & 30.4 & 2.83 & 3.47 &   22 & 1.34 \\ 
    &  $\hat \tau_{1,2|\Z=\z}^{(3)}$  &   46 & 92.8 & 42.2 & 16.5 & 42.6 & 6.45 & 10.4 & 30.4 & 2.94 & 8.06 & 22.1 &  1.4 \\ 
    &  $\hat \tau_{1,2|\Z=\z}^{(4)}$  & 15.6 &  100 & 35.2 & 5.11 &   44 & 5.31 & 3.17 &   31 & 2.45 &  3.5 & 22.4 & 1.17 \\ 
        
    \end{tabular}
    }
    \caption{Results of the simulation in Setting 2 using $h = \alpha_h \times h^{CV}$ where $h^{CV}$ has been chosen by cross-validation. All values have been multiplied by 1000. Bold values indicate optimal choices for the chosen measure of performance.}
    \label{tab:result_simulation_znorm_hCV}
    
\end{table}

\begin{table}
    \centering
    \begin{tabular}{c|cccc}
        $n$ & 100 & 500 & 1000 & 2000 \\
        \hline
        $\EE[h^{CV}]$ & 0.77 & 0.43 & 0.34 & 0.27 \\
        $Sd[h^{CV}]$ & 0.17 & 0.091 & 0.060 & 0.057 \\
        $h^{ref} = n^{-1/5}$ & 0.40 & 0.29 & 0.25 & 0.22
    \end{tabular}
    \caption{Expectation and standard deviation of the bandwidth selected by cross-validation as a function of the sample size $n$, and comparison with bandwidth $h^{ref}$ chosen by the rule-of-thumb.}
    \label{tab:summary_statistics_hCV}
\end{table}



\bibliographystyle{elsarticle-num} 
\bibliography{biblio_KendallRegression}{}


%
\appendix

\section{Proofs}

For convenience, we recall Berk's (1970) inequality (see Theorem A in Serfling \cite[p.201]{serfling1980approximation}). Note that, if $m=1$, this reduces to Bernstein's inequality.
\begin{lemma}
    Let $m,n >0$, $\X_1, \dots, \X_n$ i.i.d. random vectors with values in a measurable space $\Xc$ and $g: \Xc^m \to [a, b]$ be a symmetric real bounded function. Set $\theta := \EE[g(\X_1, \dots, \X_m)]$ and $\sigma^2 := Var[g(\X_1, \dots, \X_m)]$.
    Then, for any $t > 0$ and $n \geq m$,
    \begin{align*}
        \PP \left( \binom{n}{m}^{-1} \sum_c g(\X_{i_1}, \dots, \X_{i_m}) - \theta \geq t \right)
        \leq \exp \bigg( - \frac{[ n / m] t^2}{2 \sigma^2 + (2/3) (b-\theta) t} \bigg),
    \end{align*}
    where $\sum_c$ denotes summation over all subgroups of $m$ distinct integers $(i_1,\ldots,i_m)$ of $\{1, \dots n \}$.
    \label{lemma:bernstein_U_stat}
\end{lemma}

\subsection{Proof of Proposition \ref{prop:relationship_hat_tau_i}}
\label{proof:prop:relationship_hat_tau_i}

Since there are no ties a.s., 
\begin{align*}
    1 + \hat \tau_{1,2|\Z=\z}^{(1)}
    & = 4 \sum_{i=1}^n \sum_{j=1}^n w_{i,n}(\z) w_{j,n}(\z) \Big( 
    \1 \big\{ X_{i,1} < X_{j,1} \big\} - \1 \big\{ X_{i,1} < X_{j,1} , X_{i,2} > X_{j,2} \big\} \Big) \\
    &= 4 \sum_{i=1}^n \sum_{j=1}^n w_{i,n}(\z) w_{j,n}(\z)
    \1 \big\{ X_{i,1} < X_{j,1} \big\}+\hat \tau_{1,2|\Z=\z}^{(3)} - 1 .
\end{align*}
But 
\begin{align*}
    1 &= \sum_{i=1}^n \sum_{j=1}^n w_{i,n}(\z) w_{j,n}(\z) = \sum_{i=1}^n \sum_{j=1}^n w_{i,n}(\z) w_{j,n}(\z)  \Big(   \1 \big\{ X_{i,1} \leq X_{j,1} \big\}+  \1 \big\{ X_{i,1} > X_{j,1} \big\} \Big) \\
    &= 2 \sum_{i=1}^n \sum_{j=1}^n w_{i,n}(\z) w_{j,n}(\z)  \1 \big\{ X_{i,1} < X_{j,1} \big\} + \sum_{i=1}^n w_{i,n}^2, 
\end{align*}
implying 
$1+ \hat \tau_{1,2|\Z=\z}^{(1)}=2 (1-s_n)+\hat \tau_{1,2|\Z=\z}^{(3)} -1,$
and then $\hat \tau_{1,2|\Z=\z}^{(1)}=\hat \tau_{1,2|\Z=\z}^{(3)}-2s_n$.
\mds
Moreover, 
\begin{align*}
    \hat \tau_{1,2|\Z=\z}^{(2)}
    &=  \sum_{i=1}^n \sum_{j=1}^n w_{i,n}(\z) w_{j,n}(\z)
    \Big( \1 \big\{ X_{i,1} > X_{j,1} , X_{i,2} > X_{j,2} \big\}
    +  \1 \big\{ X_{i,1} < X_{j,1} , X_{i,2} < X_{j,2} \big\} \\
    &- \1 \big\{ X_{i,1} > X_{j,1} , X_{i,2} < X_{j,2} \big\} - \1 \big\{ X_{i,1} < X_{j,1} , X_{i,2} > X_{j,2} \big\} \\
    &= 2 \sum_{i=1}^n \sum_{j=1}^n w_{i,n}(\z) w_{j,n}(\z)
    \Big( \1 \big\{ X_{i,1} > X_{j,1} , X_{i,2} > X_{j,2} \big\}
    - \1 \big\{ X_{i,1} > X_{j,1} , X_{i,2} < X_{j,2} \big\} \Big) \\
    &= \frac{1}{2}\big(\hat \tau_{1,2|\Z=\z}^{(1)} + 1  \big)
    + \frac{1}{2}\big( \hat \tau_{1,2|\Z=\z}^{(3)} - 1  \big)
    = \frac{\hat \tau_{1,2|\Z=\z}^{(1)}+\hat \tau_{1,2|\Z=\z}^{(3)}}{2}
    =\hat \tau_{1,2|\Z=\z}^{(1)} + s_n
    = \hat \tau_{1,2|\Z=\z}^{(3)} - s_n.
    \; \Box
\end{align*}

\mds

\subsection{Proof of Proposition~\ref{cor:probaTau_Z_valid}}
\label{proof:lemma:bound_f_hat_f}

\begin{lemma}
    Under Assumptions \ref{assumpt:kernel_integral},~\ref{assumpt:f_Z_Holder} and~\ref{assumpt:f_Z_max},
    we have for any $t > 0$,
    \begin{align*}
        \PP \bigg( \big| \hat f_{\Z}(\z)-f_{\Z}(\z) \big|
        \geq \frac{  C_{K, \alpha} h^{\alpha}}{ \alpha !} + t \bigg)
        \leq 2 \exp \bigg( - \frac{n h^p t^2}{2 f_{\Z, max} \int K^2 + (2/3) C_K t} \bigg).
    \end{align*}
    \label{lemma:bound_f_hat_f}
\end{lemma}
This Lemma is proved below.
If, for some $\epsilon > 0$, we have
$  C_{K, \alpha} h^{\alpha} / \alpha ! + t
\leq f_{\Z, min} - \epsilon$,
then $\hat f(\z) \geq \epsilon > 0$ with a probability larger than
$1 - 2 \exp \big( - n h^p t^2 / (2 f_{\Z, max} \int K^2 + (2/3) C_K t) \big)$.
So, we should choose the largest $t$ as possible, which yields Proposition~\ref{cor:probaTau_Z_valid}.

\mds

It remains to prove Lemma~\ref{lemma:bound_f_hat_f}.
Use the usual decomposition between a stochastic component and a bias:
$\hat f_{\Z}(\z)-f_{\Z}(\z) = \big(\hat f_{\Z}(\z) - \EE[\hat f_{\Z}(\z)] \big)
+ \big( \EE[\hat f_{\Z}(\z)] - f_{\Z}(\z) \big)$. We first bound the bias from above.
    $$ \EE[\hat f_{\Z}(\z)] - f_{\Z}(\z)
    =  \int_{\Rb^p} K(\u) \Big(f_{\Z} \big(\z+ h \u \big) - f_{\Z}(\z) \Big) d\u.$$
Set $\phi_{\z, \u}(t) := f_{\Z} \big( \z + t h \u \big)$ for $t \in [0,1]$. This function has at least the same regularity as $f_\Z$, so it is $ \alpha $-differentiable.
By a Taylor-Lagrange expansion, we get
$$     \int_{\Rb^p} K(\u)
    \Big(f_{\Z} \big(\z+h \u \big) - f_{\Z}(\z) \Big) d\u
    =  \int_{\Rb^p} K(\u)
    \bigg(\sum_{i=1}^{ \alpha -1} \frac{1}{i!} \phi_{\z, \u}^{(i)}(0)
    + \frac{1}{ \alpha !}\phi_{\z, \u}^{( \alpha )}
    (t_{\z, \u}) \bigg) d\u ,$$
for some real number $t_{\z, \u}\in (0,1)$.
By Assumption \ref{assumpt:kernel_integral} and for every
$i<\alpha$, $\int_{\Rb^p} K(\u) \phi_{\z, \u}^{(i)}(0) \,d\u = 0$. Therefore,
\begingroup \allowdisplaybreaks
\begin{eqnarray*}
    \lefteqn{
    \Big| \EE[\hat f_{\Z}(\z)] - f_{\Z}(\z) \Big| = \Big| \int_{\Rb^p} K(\u) \frac{1}{ \alpha !}
    \phi_{\z, \u}^{( \alpha )}(t_{\z, \u}) d\u \Big| }\\
    &=& \frac{1}{\alpha  !}
    \Big| \int_{\Rb^p} K(\u) \sum_{i_1, \dots, i_{ \alpha } = 1}^{p}
    h^{ \alpha } u_{i_1} \dots u_{i_{\alpha }}
    \frac{\partial^{\alpha } f_{\Z}}
    {\partial z_{i_1} \dots  \partial z_{i_{ \alpha }}}
    \big( \z + t_{\z, \u} h \u \big) d\u \Big|
    \leq \frac{ C_{K, \alpha}}{ \alpha  !} h^{\alpha}.
\end{eqnarray*}
\endgroup

\medskip

Second, the stochastic component may be written as
\begin{align*}
     \hat f_{\Z}(\z) - \EE[\hat f_{\Z}(\z)]
    &=  n^{-1} \sum_{i=1}^n K_h(\Z_i-\z)
    - \EE \Big[n^{-1} \sum_{i=1}^n K_h(\Z_i-\z) \Big]
    =  n^{-1} \sum_{i=1}^n \big( g_\z(\Z_i) - \EE[g(\Z_i)]\big),
\end{align*}
where $g(\Z_i) := K_h(\Z_i-\z)$.
Apply Lemma \ref{lemma:bernstein_U_stat} with $m = 1$ and the latter $g(\Z_i)$.
Here, we have $ b = -a = h^{-p} C_K$, $\theta = \EE[g(\Z_1)]\geq 0 $ and $\big| Var[g(\Z_1)] \big| \leq h^{-p} f_{\Z, max} \int K^2$, and we get
\begin{align*}
    \PP \left( \big| \frac{1}{n} \sum_{i=1}^n K_h(\Z_i-\z) - \EE[K_h(\Z_i-\z)] \big| \geq t \right) \leq 2\exp \bigg( - \frac{n t^2}{2 h^{-p} f_{\Z, max} \int K^2 + (2/3) h^{-p} C_K t} \bigg).
    \; \Box
\end{align*}

\mds

\subsection{Proof of Proposition \ref{prop:exponential_bound_KendallsTau}}
\label{proof:prop:exponential_bound_KendallsTau}
We show the result for $k=1$. The two other cases can be proven in the same way. 

Consider the decomposition
\begingroup \allowdisplaybreaks
\begin{eqnarray*}
\lefteqn{      \hat \tau_{1,2|\Z=\z} - \tau_{1,2|\Z=\z}     =  4\sum_{1 \leq i,j \leq n} w_{i,n}(\z) w_{j,n}(\z)
    \1 \big\{ \X_{i} < \X_{j} \big\} - 4\PP \big( \X_{1} < \X_{2} \big| \Z_1 = \Z_2 = \z \big)  }\\
    &=& \frac{4}{n^2\hat f_{\Z}^2(\z)}   \sum_{1 \leq i,j \leq n}
    K_h(\Z_i-\z) K_h(\Z_j-\z)
    \Big( \1 \big\{ \X_{i} < \X_{j} \big\}
    - \PP \big( \X_{1} < \X_{2} \big| \Z_1 = \Z_2 = \z \big) \Big) \\
    &=:& \frac{4}{\hat f_{\Z}^2(\z)} \sum_{1 \leq i,j \leq n} S_{i,j}(\z).
\end{eqnarray*}
\endgroup
Therefore, for any positive numbers $x$ and $\lambda(\z)$, we have
\begin{eqnarray*}
\lefteqn{   \PP(    |\hat \tau_{1,2|\Z=\z} - \tau_{1,2|\Z=\z} | > x) \leq 
    \PP\Big(    \frac{1}{\hat f_{\Z}^2(\z)}   > \frac{1+\lambda(\z)}{f_{\Z}^2(\z)} \Big) 
    +      \PP\Big(    \frac{4(1+\lambda(\z))}{f_{\Z}^2(\z)}   \times | \sum_{1 \leq i,j \leq n} S_{i,j}(\z) | > x  \Big) }\\
&\leq & \PP\Big(   | \frac{1}{\hat f_{\Z}^2(\z)} - \frac{1}{f_{\Z}^2(\z)}  |  > \frac{\lambda(\z)}{f_{\Z}^2(\z)} \Big) 
+      \PP\Big(    \frac{4(1+\lambda(\z))}{f_{\Z}^2(\z)}   \times | \sum_{1 \leq i,j \leq n} S_{i,j}(\z) | > x  \Big).\hspace{5cm}
\end{eqnarray*}
For any $t$ s.t. $  C_{K, \alpha} h^{\alpha} /  \alpha  ! + t < f_{\Z, min}/2$, set
$\lambda(\z)=16 f^2_{\z}(\z) \big(  C_{K, \alpha} h^{\alpha}/ \alpha  ! + t \big)/f_{\Z, min}^3.$
This yields
\begin{eqnarray*}
\lefteqn{   \PP\Big(    |\hat \tau_{1,2|\Z=\z} - \tau_{1,2|\Z=\z} | > x\Big) \leq 
     \PP\Big(   | \frac{1}{\hat f_{\Z}^2(\z)} - \frac{1}{f_{\Z}^2(\z)}  |  > \frac{16} {f_{\Z, min}^3} \big( \frac{  C_{K, \alpha} h^{\alpha}} { \alpha  !} + t \big) \Big) }\\
&+&       \PP\Big(    | \sum_{1 \leq i,j \leq n} S_{i,j}(\z) | >  \frac{f_{\z}^2(\z) x}{4(1+\lambda(\z))}  \Big).\hspace{8cm}
\end{eqnarray*}
By setting 
$$ x=\frac{4}{f^2_{\z}(\z)}\Big( \frac{C_{\X\Z, \alpha}   h^\alpha}{ \alpha  !}+ \frac{3f_{\z}(\z)\int K^2}{2nh^p} + t' \Big)
\bigg( 1+ \frac{16 f^2_{\Z}(\z)} {f_{\Z, min}^3} 
\Big( \frac{  C_{K, \alpha} h^{\alpha}} { \alpha  !} + t \Big) \bigg),$$
and applying the next two lemmas~\ref{lemma:bound_f2_hatf2} and~\ref{exp_bound_Sijz}, we get the result.
$\;\; \Box$

\begin{lemma}
    Under Assumptions \ref{assumpt:kernel_integral}-\ref{assumpt:f_Z_max}
    and if $  C_{K, \alpha} h^{\alpha} /  \alpha  ! + t < f_{\Z, min}/2$
    for some $t > 0$,
    \begin{align*}
        \PP \bigg( | \frac{1}{\hat f_{\Z}^2(\z)} - \frac{1}{ f_{\Z}^2(\z)}  | >
        \frac{16} {f_{\Z, min}^3}
        \Big( \frac{  C_{K, \alpha} h^{\alpha}}
        { \alpha  !} + t \Big) \bigg)
        \leq 2 \exp \bigg( - \frac{n h^p t^2}{2 f_{\Z, max} \int K^2 + (2/3) C_K t} \bigg),
    \end{align*}
    and $\hat f_{\Z}(\z)$ is strictly positive on these events.
    \label{lemma:bound_f2_hatf2}
\end{lemma}

{\it Proof :} Applying the mean value inequality to the function $x \mapsto 1/x^2$, we get
the inequality
$\Big|  1/\hat f_{\Z}^2(\z)  -  1/{f_{\Z}^2(\z)}  \Big|
\leq  2
\big| \hat f_{\Z}(\z)-f_{\Z}(\z) \big|/ f_{\Z}^*{}^3,$
where $f_{\Z}^*$ lies between $\hat f_{\Z}(\z)$ and $f_{\Z}(\z)$.
Denote by $\Ec$ the event 
$\Ec:=\big\{ | \hat f_{\Z}(\z)-f_{\Z}(\z) |
    \leq   C_{K, \alpha} h^{\alpha} /\alpha  ! + t \big\}$.
By Lemma \ref{lemma:bound_f_hat_f}, we obtain
\begin{align}
    \PP ( \Ec )
    \geq 1 - 2 \exp \Big( - \frac{n h^p t^2}{2 f_{\Z, max} \int K^2 + (2/3) C_K t} \Big).
    \label{arg_exp_bound_f2_hatf2}
\end{align}
Therefore, on this event $\Ec$, $\big| \hat f_{\Z}(\z)-f_{\Z}(\z) \big| \leq f_{\Z, min}/2$,
so that $f_{\Z, min}/2 \leq \hat f_{\Z}(\z)$. We have also $f_{\Z, min}/2 \leq f_{\Z}(\z)$ and then $f_{\Z, min}/2 \leq f_{\Z}^*$.
Combining the previous inequalities, we finally get
\begin{align*}
    \bigg| \frac{1}{\hat f_{\Z}^2(\z)}-\frac{1}{f_{\Z}^2(\z)} \bigg|
    &\leq \frac{16}{f_{\Z, min}^3} \big| \hat f_{\Z}(\z)-f_{\Z}(\z) \big|
    \leq \frac{16}{f_{\Z, min}^3}
    \Big( \frac{  C_{K, \alpha} h^{\alpha}} { \alpha  !} + t \Big),
\end{align*}
on $\Ec$. But since
    $$     \PP \bigg( | \frac{1}{\hat f_{\Z}^2(\z)} - \frac{1}{f_{\Z}^2(\z)} | > 
        \frac{16 }{f_{\Z, min}^3} \Big( \frac{  C_{K, \alpha} h^{\alpha}}
        { \alpha  !} + t \Big) \bigg)  \leq  \PP(\Ec^c) ,$$
we deduce the result.$\;\; \Box$
\medskip

\begin{lemma}
    Under Assumptions \ref{assumpt:kernel_integral}-\ref {assumpt:f_XZ_Holder}, if $C_{\tilde K,2} h^{2} < f_{\z}(\z)$, we have for any $t>0$ 
    \begin{eqnarray*}
 \lefteqn{        \PP \bigg( \Big| \sum_{1 \leq i,j \leq n} S_{i,j}(\z) \Big| >
        \frac{C_{\X\Z, \alpha}   h^\alpha}
        { \alpha  !} + \frac{3f_{\z}(\z)\int K^2}{2nh^p}+ t \bigg)
        \leq  2 \exp \bigg( - \frac{(n-1) h^{2p} t^2 }{4 f_{\Z, max}^2 (\int K^2)^2 + (8/3) C_K^2  t} \bigg) }\\
        &+& 2 \exp \bigg( - \frac{n h^p (f_{\z}(\z) - C_{\tilde K,2}h^2)^2 }{8 f_{\Z, max} \int \tilde K^2 + 4C_{\tilde K} (f_{\z}(\z) - C_{\tilde K,2}h^2)/3} \bigg).
        \hspace{8cm}
    \end{eqnarray*}
\label{exp_bound_Sijz}
\end{lemma}
{\it Proof :}
Note that
$\sum_{1 \leq i,j \leq n} S_{i,j}(\z)
= \sum_{1 \leq i\neq j \leq n} \big( S_{i,j}(\z) - \EE[S_{i,j}(\z)] \big)
+ n(n-1)\EE[S_{1,2}(\z)] + \sum_{i=1}^n S_{i,i}(\z).$
The ``diagonal term''
$\sum_{i=1}^n S_{i,i}(\z) = -  \PP \big( \X_1 < \X_2 \big| \Z_1 = \Z_2 = \z \big)
\sum_{i=1}^n K^2_h(\Z_i-\z) / n^2 $ is negative and negligible. It will be denoted by $-\Delta_n(\z) <0$.
Note that $\tilde K(\cdot):= K^2(\cdot)/\int K^2$ is a two-order kernel. Then, $\tilde f_{\z}(\z):= \sum_{i=1}^n \tilde K_h(\Z_i-\z) / n $ is a consistent estimator of $f_{\Z}(\z)$. Therefore, due to Lemma~\ref{lemma:bound_f_hat_f} and with obvious notations, we have for every $\eps>0$
    \begin{align*}
        \PP \bigg( \big| \tilde f_{\Z}(\z)-f_{\Z}(\z) \big|
        \geq \frac{  C_{\tilde K,2} h^{2}}{2} + \eps \bigg)
        \leq 2 \exp \bigg( - \frac{n h^p \eps^2 }{2 f_{\Z, max} \int \tilde K^2 + (2/3) C_{\tilde K} \eps} \bigg).
    \end{align*}
This implies 
    \begin{eqnarray*}
    \lefteqn{    \PP \bigg( | \frac{\int K^2}   {n^2 h^p } \sum_{i=1}^n \tilde K_h(\Z_i-\z) - \frac{ f_{\Z}(\z) \int K^2}{n h^p} | \geq 
        \Big(\frac{\int  K^2}{ n h^p} \Big) \Big(\frac{  C_{\tilde K,2} h^{2}}{2} + \eps \Big) \bigg)   }\\
        &\leq  & 2 \exp \bigg( - \frac{n h^p \eps^2 }{2 f_{\Z, max} \int \tilde K^2 + (2/3) C_{\tilde K} \eps} \bigg). \hspace{6cm}
    \end{eqnarray*}
By choosing $\eps$ s.t. $C_{\tilde K,2} h^{2}/2 + \eps = f_{\z}(\z)/2$, $\Delta_n$ will be smaller than $3f_{\z}(\z) \int K^2 /(2nh^p)$ with a probability that is larger than 
\begin{equation}
1- 2 \exp \bigg( - \frac{n h^p \eps^2 }{2 f_{\Z, max} \int \tilde K^2 + (2/3) C_{\tilde K} \eps} \bigg).
\label{proba_delta_n}    
\end{equation}

\mds

Now, let us deal with the main term, that is decomposed as a stochastic component and a bias component.
First, let us deal with the bias. Simple calculations provide, if $i\neq j$,
\begingroup \allowdisplaybreaks
\begin{align*}
    &\EE[S_{i,j}(\z)]
    = n^{-2}\EE \bigg[ K_h(\Z_i-\z) K_h(\Z_j-\z)
    \Big( \1 \big\{ \X_{i} < \X_{j} \big\}
    - \PP \big( \X_{i} < \X_{j} \big| \Z_i = \Z_j = \z \big) \Big) \bigg] \\
    &= n^{-2}\int_{\Rb^{2p+2}} K_h(\z_1-\z) K_h(\z_2-\z)
    \Big( \1 \big\{ \x_1 < \x_2 \big\}
    - \PP \big( \X_{i} < \X_{j} \big| \Z_i = \Z_j = \z \big) \Big) \\
    & \hspace{2cm} \times  f_{\X,\Z}(\x_1,\z_1) \, f_{\X, \Z}(\x_2,\z_2) \,
    d\x_1 \,d\z_1\, d\x_2 \, d\z_2 \\
    &= n^{-2}\int_{\Rb^{2p+2}} K(\u) K(\v)
    \Big( \1 \big\{ \x_1 < \x_2 \big\}
    - \PP \big( \X_{i} < \X_{j} \big| \Z_i = \Z_j = \z \big) \Big) \\
    & \hspace{2cm} \times   \bigg( f_{\X,\Z} \Big(\x_1,\z + h \u \Big) \,
    f_{\X,\Z} \Big(\x_2,\z + h \v \Big)
    - f_{\X,\Z} (\x_1,\z) \, f_{\X,\Z} (\x_2,\z) \bigg)\, d\x_1\, d\u\, d\x_2\, d\v,
\end{align*}
\endgroup
because, for every $\z$,
$$    0
    = \int_{\Rb^4} \Big( \1 \big\{ \x_1 < \x_2 \big\}
    - \PP \big( \X_1 < \X_2 \big| \Z_1 = \Z_2 = \z \big) \Big)
    \, f_{\X,\Z}(\x_1,\z) f_{\X, \Z}(\x_2,\z)\, d\x_1 \, d\x_2. $$
Apply the Taylor-Lagrange formula to the function
$\phi_{\x_1, \x_2, \u, \v}(t) := f_{\X,\Z} \big(\x_1,\z + t h \u \big) \,
f_{\X,\Z} \big(\x_2,\z + t h \v \big).$
With obvious notation, this yields
\begingroup \allowdisplaybreaks
\begin{align*}
    \EE[S_{i,j}(\z)]
    &= n^{-2}\int K(\u) K(\v)
    \Big( \1 \big\{ \x_1 < \x_2 \big\}
    - \PP \big( \X_{i} < \X_{j} \big| \Z_i = \Z_j = \z \big) \Big) \\
    & \hspace{2cm} \times \bigg( \sum_{k=1}^{ \alpha - 1} \frac{1}{k !}
    \phi_{\x_1, \x_2, \u, \v}^{(k)}(0) + \frac{1}{ \alpha  !}
    \phi_{\x_1, \x_2, \u, \v}^{(\alpha )}
    (t_{\x_1, \x_2, \u, \v}) \bigg) \, d\x_1\, d\u\, d\x_2\, d\v \\
    &= \int \frac{K(\u) K(\v)}{n^2 \alpha  ! }
    \Big( \1 \big\{ \x_1 < \x_2 \big\}
    - \PP \big( \X_{i} < \X_{j} \big| \Z_i = \Z_j = \z \big) \Big)
    \phi_{\x_1, \x_2, \u, \v}^{( \alpha )}
    (t_{\x_1, \x_2, \u, \v})  d\x_1 \, d\u\, d\x_2 \,d\v .
\end{align*}
\endgroup
Since $\phi_{\x_1, \x_2, \u, \v}^{( \alpha )} (t) $ is equal to
$$\sum_{k=0}^{ \alpha } \binom{ \alpha }{k}
    \sum_{i_1, \dots, i_{ \alpha } = 1}^{p}
    h^{ \alpha } u_{i_1} \dots u_{i_k}
    v_{i_{k+1}} \dots v_{i_{ \alpha }}
    \frac{\partial^{k} f_{\X, \Z}}
    {\partial z_{i_1} \dots  \partial z_{i_k}} \Big( \x_1, \z + t h \u \Big)
    \frac{\partial^{ \alpha -k} f_{\X, \Z}}
    {\partial z_{i_{k+1}} \dots  \partial z_{i_{ \alpha }}}
    \Big( \x_2, \z + t h \v \Big) ,$$
using Assumption \ref{assumpt:f_XZ_Holder}, we get
\begin{equation}
\big| \EE[S_{1,2}(\z)] \big| \leq C_{\X\Z, \alpha}   h^\alpha
    /(n^2  \alpha  !).
\label{upper_bound_ES12}
\end{equation}

\medskip

Second, the stochastic component will be bounded from above. Indeed,
\begin{align*}
    & \sum_{1 \leq i\neq j \leq n} (S_{i,j}(\z) - \EE[S_{i,j}(\z)])
    = \frac{1}{n^2}
     \sum_{1 \leq i \neq j \leq n}
    g_\z \big((\X_i, \Z_i) \, , \, (\X_j, \Z_j) \big) ,
\end{align*}
with the function $g_\z$ defined by
\begin{eqnarray*}
    \lefteqn{ g_\z \big((\X_i, \Z_i), (\X_j, \Z_j) \big)
    :=  K_h(\Z_i-\z) K_h(\Z_j-\z)
    \Big( \1 \big\{ \X_{i} < \X_{j} \big\}
    - \PP \big( \X_{i} < \X_{j} \big| \Z_i = \Z_j = \z \big) \Big)  }\\
    & - \EE \bigg[ K_h(\Z_i-\z) K_h(\Z_j-\z)
    \Big( \1 \big\{ \X_{i} < \X_{j} \big\}
    - \PP \big( \X_{i} < \X_{j} \big| \Z_i = \Z_j = \z \big) \Big) \bigg].\hspace{2cm}
\end{eqnarray*}
%
The symmetrized version of $g$ is
$\tilde g_{i,j}=\Big(
g_\z \big((\X_i, \Z_i) \, , \, (\X_j, \Z_j) \big)
+ g_\z \big((\X_j, \Z_j) \, , \, (\X_i, \Z_i) \big)  \Big)/2.$
We can now apply Lemma \ref{lemma:bernstein_U_stat} to the sum of the $\tilde g_{i,j}$.
With its notation, $\theta = \EE \big[\tilde g_{i,j}  \big] = 0$. Moreover,
\begin{eqnarray*}
\lefteqn{  \Big| Var \Big[g_\z \big((\X_i, \Z_i),(\X_j, \Z_j) \big) \Big] \Big| }\\
& \leq & \int K^2_h(\z_1-\z) K^2_h(\z_2-\z)
    \Big( \1 \big\{ \x_{1} < \x_{2} \big\}
    - \PP \big( \X_{i} < \X_{j} \big| \Z_i = \Z_j = \z \big) \Big)^2 \\
    &\times & f_{\X,\Z}(\x_1,\z_1) f_{\X,\Z}(\x_2,\z_2) \, d\x_1\, d\x_2 \, d\z_1\, d\z_2 \\
&\leq &
\int \frac{K^2(\t_1) K^2(\t_2)}{h^{2p}} f_{\X,\Z}(\x_1,\z -h\t_1) f_{\X,\Z}(\x_2,\z - h\t_2) \, d\x_1\, d\x_2
    \, d\t_1\, d\t_2 \\
&\leq & h^{-2p} f_{\Z, max}^2 \Big(\int K^2\Big)^2,
\end{eqnarray*}
and the same upper bound applies for $\tilde g_{i,j}$ (invoke Cauchy-Schwarz inequality).
Here, we choose $b = -a = 2 C_K^2 h^{-2p} $.
This yields
\begin{align}
    \PP \Big( \frac{2}{n(n-1)} \sum_{1 \leq i < j \leq n} \tilde g_{i,j} > t\Big)
    \leq \exp \bigg( - \frac{[n / 2] t^2 }{2 h^{-2p} f_{\Z, max}^2 (\int K^2)^2 + (4/3) C_K^2 h^{-2p}  t} \bigg)
    \label{eq:bound_bernestein_improved}
\end{align}
Then, for every $t > 0$, we obtain
\begin{eqnarray*}
    \lefteqn{  \PP \Big( | \sum_{1 \leq i \neq j \leq n}
    \big( S_{i,j}(\z) - \EE[S_{i,j}(\z)] \big) |
    \geq t \Big) \leq  \PP \Big( \frac{1}{n^2}
    | \sum_{1 \leq i \neq j \leq n}
    g_\z \big((\X_i, \Z_i) \, , \, (\X_j, \Z_j)  | \big)
    \geq t  \Big) }\\
    & \leq &\PP \Big( \frac{(n-1)}{n} \times \frac{2}{n(n-1)}
     | \sum_{1 \leq i < j \leq n}
    \tilde g_{i,j} |
    \geq t  \Big) 
     \leq  2 \exp \bigg( - \frac{[n/ 2] t^2 }{2 h^{-2p} f_{\Z, max}^2 (\int K^2)^2 + (4/3) C_K^2 h^{-2p}  t} \bigg). 
\end{eqnarray*}
The latter inequality,~(\ref{proba_delta_n}) and~(\ref{upper_bound_ES12}) yield the result. $\;\;\Box$

\mds

\subsection{Proof of Proposition~\ref{prop:exponential_bound_KendallsTau_major}}
\label{proof:prop:exponential_bound_KendallsTau_major}
With the notations of the proof of Proposition~\ref{prop:exponential_bound_KendallsTau}, we get the following lemma, that straightforwardly implies the result.
\begin{lemma}
    Under the Assumptions and conditions of Proposition~\ref{prop:exponential_bound_KendallsTau_major},
    we have
    \begin{eqnarray*}
    \lefteqn{     \PP \bigg( \Big| \sum_{1 \leq i,j \leq n} S_{i,j}(\z) \Big| >
        \frac{C_{\X\Z, \alpha}   h^\alpha}
        { \alpha  !}+ \frac{3f_{\z}(\z)\int K^2}{2nh^p} + t \bigg)   
        \leq  C_2\exp \bigg( - \frac{\alpha_2 n h^{p} t}{  8 f_{\Z, max} (\int K^2)} \bigg)  }\\
        &+&   2 \exp \bigg( - \frac{n h^p (f_{\z}(\z) - C_{\tilde K,2}h^2)^2 }
        {8 f_{\Z, max} \int \tilde K^2 + 4 C_{\tilde K} (f_{\z}(\z) - C_{\tilde K,2}h^2)/3} \bigg) \\
        &+&
        2 \exp \Big( \frac{nh^p t^2}{32  \int K^2 (\int |K|)^2 f^3_{\Z,max} + 8 C_K\int |K|f_{\Z,max} t/3}    \Big). \hspace{7cm}
    \end{eqnarray*}
\label{exp_bound_Sijz_major}
\end{lemma}

{\it Proof :} We lead exactly the same reasoning and the same notations as in Lemma~\ref{exp_bound_Sijz}, until~(\ref{eq:bound_bernestein_improved}).
Now, with the same notations, introduce $\bar g_i:=\EE[\tilde g_{i,j} | \X_i,\Z_i]$ and consider $ \xi_{i,j}:= \tilde g_{i,j} - \bar g_i - \bar g_j.$ Then, $\xi_{i,j}$ is a degenerate (symmetrical) U-statistics because $ \EE[\xi_{i,j} | \X_i,\Z_i]=\EE[\xi_{i,j} | \X_j,\Z_j]=0, $ when $i\neq j$.
Actually $\xi_{i,j} =:\xi_\z(\X_i,\Z_i,\X_j,\Z_j)$ for some function $\xi_\z$
and set
\begin{equation}
 \ell_{\z}:(\x_1,\z_1,\x_2,\z_2)\mapsto  \frac{h^{2p}}{4C_K^2}\xi_\z \big((\x_1, \z_1) \, , \, (\x_2, \z_2) \big),
 \label{def_ell_z}
 \end{equation}
 for a fixed $\z$ and a fixed $h$.
This yields $\| \ell_{\z} \|_{\infty} \leq 1$. By usual changes of variables, we get
\begin{eqnarray*} 
\lefteqn{ 
\int \ell_{\z}^2(\x_1,\z_1,\x_2,\z_2) \, f_{\X,\Z}(\x_1,\z_1) f_{\X,\Z}(\x_2,\z_2) \, d\x_1\, d\x_2 \, d\z_1\, d\z_2 }\\
&\leq & 3h^{2p}\frac{ (\int K^2 f_{\z,max})^2}{(4C_K^2)^2}  
+  6 h^{3p}\frac{\int K^2 f_{\z,max} (\int |K| f_{\z,max})^2}{(4C_K^2)^2 } \leq \sigma^2,
\;\text{with}
\end{eqnarray*}
\begin{equation}
 \sigma := h^p C_\sigma, \; C_\sigma:= \int K^2 f_{\z,max}/ ( 2 C_K^2),
 \label{value_sigma}
 \end{equation}
because $ 6 h^{p}\int K^2 f_{\z,max} (\int |K| f_{\z,max})^2
\leq  (\int K^2 f_{\z,max})^2$.
With the notations of~\cite{major2006estimate}, this implies $D=1$, $m=1$ and $L$ is arbitrarily small. Therefore, Theorem 2 in~\cite{major2006estimate} yields
\begin{align}
    \PP \Big( \frac{1}{2n}| \sum_{i\neq j} \ell_{\z}(\X_i,\Z_i,\X_j,\Z_j) |   > x\Big)
    \leq C_2\exp \bigg( - \frac{ \alpha_2 x}{ \sigma} \bigg),
    \label{eq:bound_major0}
\end{align}
for some universal constants $C_2$ and $\alpha_2$ when $x \leq n \sigma^3$. 
By setting $t/2= 4 C_K^2 x/ (nh^{2p}) $ and applying Lemma~\ref{lemma:bernstein_U_stat}, this provides
\begin{eqnarray*}
\lefteqn{ \PP \Big( | \sum_{1 \leq i \neq j \leq n}
    \big( S_{i,j}(\z) - \EE[S_{i,j}(\z)]\big) |
    \geq t \Big) \leq 
    \PP \Big( \frac{1}{n^2}| \sum_{1 \leq i \neq j \leq n} \xi_{ij}  |
    \geq t/2 \Big) +
    \PP \Big( | \frac{1}{n} \sum_{i=1}^n \bar g_i  |
    \geq t/4 \Big)   }\\
    &\leq &
    C_2\exp \Big( - \frac{\alpha_2 n t h^{p}}{ 8  f_{\Z, max} (\int K^2)} \Big) 
    +  2 \exp \Big( \frac{nh^p t^2}{32  \int K^2 (\int |K|)^2 f^3_{\Z,max} + 8/3 C_K\int |K|f_{\Z,max} t}    \Big),
\end{eqnarray*}
when $ t \leq 2 h^p  (\int K^2)^3 f^3_{\Z,max} /   C_K^4$.
This inequality,~(\ref{proba_delta_n}) and~(\ref{upper_bound_ES12}) conclude the proof. $\;\Box$

\mds

\subsection{Proof of Proposition~\ref{prop:exponential_bound_KendallsTau_uniform}}
\label{proof:prop:exponential_bound_KendallsTau_uniform}
For $k=1$, we follow the paths of the proof of Proposition~\ref{prop:exponential_bound_KendallsTau_major}.
Since 
$\hat \tau_{1,2|\Z=\z} - \tau_{1,2|\Z=\z} = 4 \sum_{1 \leq i,j \leq n} S_{i,j}(\z)/\hat f_{\Z}^2(\z)$,
we prove the result if we bound from above  
$1 / \hat f_{\Z}^2(\z)$ and $\big|\sum_{1 \leq i,j \leq n} S_{i,j}(\z) \big|$ uniformly w.r.t. $\z\in \Zc$. 
To be specific, for any positive constant $\mu <1$, if 
$| \hat f_{\Z}(\z) - f_{\Z}(\z)|\leq \mu f_{\z,min} $, then $1/\hat f^2_{\z,max}(\z)\leq f^{-2}_{\z,min}(1-\mu)^{-2}$. We deduce
\begin{eqnarray*}
\lefteqn{   \PP(\sup_{\z\in\Zc}    |\hat \tau_{1,2|\Z=\z} - \tau_{1,2|\Z=\z} | > x) \leq 
 \PP\big( \| \hat f_{\Z} - f_{\Z} \|_{\infty} > \mu f_{\z,min}  \big) }\\
 &+&      \PP(    \frac{4}{f_{\Z,min}^2(1-\mu)^2}  \sup_{\z\in\Zc}  | \sum_{1 \leq i,j \leq n} S_{i,j}(\z) | > x  ).
\end{eqnarray*}
First invoke the uniform exponential inequality, as stated in~\cite{rinaldo2010generalized}, Proposition 9: for every $\eps< b_K \int K^2 f_{\Z,max}/C_{K}$,
\begin{equation}
    \PP \big( \| \hat f_{\Z}-f_{\Z} \|_{\infty}
    > \eps + \frac{C_{\X\Z, \alpha}   h^\alpha}{\alpha  !}  \big) \leq \PP \big( \| \hat f_{\Z}-\EE[\hat f_{\Z} ] \|_{\infty} > \eps \big)
    \leq  L_K \exp \big( - C_{f,K} nh^p \eps^2 \big),
\label{exp_ineg_kernel_desity_unif}
\end{equation}
for $n$ sufficiently large. Then, apply Lemma~\ref{exp_bound_Sijz_uniform}, by setting $(x,\eps)$ so that
$$x=\frac{4}{f^2_{\z,min}(1-\mu)^2}\Big( \frac{C_{\X\Z, \alpha}   h^\alpha}{ \alpha  !} + \frac{3f_{\z,max}\int K^2}{2nh^p} + t\Big)\;\text{and}\;
\eps + \frac{ C_{\X\Z, \alpha}   h^\alpha}{\alpha  !}  =  \mu f_{\z,min}.\;\; \Box $$

\begin{lemma}
    Under the assumptions and conditions of Proposition~\ref{prop:exponential_bound_KendallsTau_uniform}, we have
    \begin{eqnarray*}
    \lefteqn{     \PP \bigg( \sup_{\z\in\Zc}\Big| \sum_{1 \leq i,j \leq n} S_{i,j}(\z) \Big| >
        \frac{C_{\X\Z, \alpha}   h^\alpha}
        { \alpha  !} + \frac{3f_{\z,max}\int K^2}{2nh^p} + t \bigg)   }\\
        &\leq &     C_2 D \exp \bigg( - \frac{\alpha_2 n t h^{p}}{  8 f_{\Z, max} (\int K^2)} \bigg) 
        +   L_{\tilde K} \exp \big( - C_{f, \tilde K} nh^p (f_{\z,max} - \tilde C_{\X\Z, 2}  h^{2})^2/4 \big) \\
        &+&
        2 \exp\big( -\frac{A_2 n h^p t^2 C_K^{-4}}{16^2 A_1^2  \int K^2 f_{\z,max}^3 (\int |K|)^2}  \big) + 2\exp\big(-\frac{A_2 n h^p t}{16 C_K^2 A_1} \big) .
    \end{eqnarray*}
\label{exp_bound_Sijz_uniform}
\end{lemma}
{\it Proof :}
We will use the arguments and notations of the proof of Lemmas~\ref{exp_bound_Sijz} and~\ref{exp_bound_Sijz_major}.
We still invoke the decomposition
$\sum_{1 \leq i,j \leq n} S_{i,j}(\z)
= \sum_{1 \leq i\neq j \leq n} \big( S_{i,j}(\z) - \EE[S_{i,j}(\z)] \big)
+ n(n-1)\EE[S_{1,2}(\z)] + \sum_{i=1}^n S_{i,i}(\z).$
First let us find a uniform bound for the ``diagonal term'' $\Delta_n(\z)=\sum_{i=1}^n S_{i,i}(\z)=\int K^2 \tilde f_{\z}(\z) /(n h^p)$.
As in~(\ref{exp_ineg_kernel_desity_unif}), for every $\eps< b_{\tilde K} \int \tilde K^2 f_{\Z,max}/C_{\tilde K}$,
\begin{align*}
    \PP \big( \| \tilde f_{\Z}-f_{\Z} \|_{\infty}
    > \eps + \frac{\tilde C_{\X\Z, 2}   h^2}{2 }\big)
    \leq  L_{\tilde K} \exp \big( - C_{f, \tilde K} nh^p \eps^2 \big),
\end{align*}
for $n$ sufficiently large. This implies
$$\PP \bigg( \sup_{\z\in\Zc}| \frac{\int K^2}{n^2 h^p } \sum_{i=1}^n \tilde K_h(\Z_i-\z) - \frac{ f_{\Z}(\z) \int K^2}{n h^p} | \geq 
\Big(\frac{\int  K^2}{ n h^p} \Big)  \big( \eps + \frac{\tilde C_{\X\Z, 2}   h^2}{2 }\big) \bigg) 
\leq    L_{\tilde K} \exp \big( - C_{f, \tilde K} nh^p \eps^2 \big).$$

Choose $\eps$ s.t. $\tilde C_{\X\Z, 2} h^{2}/2 + \eps = f_{\z,max}/2$. Then, $\sup_{\z}|\Delta_n(\z)|$ will be smaller than $3f_{\z,max} \int K^2 /(2nh^p)$ with a probability that is larger than 
\begin{equation}
1- L_{\tilde K} \exp \big( - C_{f, \tilde K} nh^p \eps^2 \big).
\label{proba_delta_n_uniform}    
\end{equation} 

Moreover, it is easy to see that
\begin{equation}
\sup_{\z\in\Zc}\big| \EE[S_{1,2}(\z)] \big| \leq C_{\X\Z, \alpha}   h^\alpha
    /(n^2   \alpha  !).
\label{upper_bound_ES12_uniform}
\end{equation}

With the notations of Lemma~\ref{exp_bound_Sijz_major}'s proof, the stochastic component is driven by
\begin{align*}
    & \sum_{1 \leq i\neq j \leq n} (S_{i,j}(\z) - \EE[S_{i,j}(\z)])
    = \frac{1}{n^2}
     \sum_{1 \leq i \neq j \leq n}
    g_\z \big((\X_i, \Z_i) \, , \, (\X_j, \Z_j) \big) \\
    &= \frac{1}{n^2} \sum_{1 \leq i \neq j \leq n} \tilde g_{i,j}
    = \frac{1}{n^2} \sum_{1 \leq i \neq j \leq n} \xi_{i,j} + 
    \frac{2(n-1)}{n^2}\sum_{i=1}^n \bar g_i  .
\end{align*}

Now apply Theorem 1 in ~\cite{major2006estimate}, by recalling~(\ref{def_ell_z}) and considering the family 
$ \Fc := \Big\{\ell_{\z},\, \z\in \Zc\Big\},$ for a fixed bandwidth $h$. The constant $\sigma$ has the same value as in~(\ref{value_sigma}).
It is easy to check that the latter class of functions is $L^2$ dense (see~\cite{major2006estimate}).
Set $\eps \in (0,1)$. Since $K$ is $\lambda_K$-Lipschitz, every function $\ell_{\z}\in\Fc$ can be approximated in $L^2$ by a function $\ell_{\z_j}\in \Fc$, for some $j \in\{1,\ldots,m\}$ s.t. $\int | \ell_{\z} - \ell_{\z_i} |^2 d\nu\leq  \eps^{2}$, for any probability measure $\nu$. 
Indeed, $ \int | \ell_{\z} - \ell_{\z_i} |^2 d\nu \leq 64\lambda^2_K \| \z-\z_j\|^2_\infty C_K^2 h^{-2}$ that is less than $\eps^{2}$, if we cover $\Zc$ by a grid of $m$ points $(\z_j)$ in $\Zc$ s.t. $\|\z-\z_j\|_\infty\leq \eps h /(8 C_K \lambda_K):=\eps\delta$.  
This can be done with $m \leq \eps^{-p} \lceil \prod_{k=1}^p \big((b_k-a_k)/\delta \big) \rceil= \eps^{-p}\lceil  \Vc \delta^{-p} \rceil$ points.
Then, with the notations of~\cite{major2006estimate}, $L=p$ and $D=\Vc (8 C_K \lambda_K/h)^{p} $.
As above, this yields
\begin{align}
    \PP \Big( \sup_{\z\in\Zc} \frac{1}{n^2} | \sum_{1 \leq i \neq j \leq n} \xi_{\Z}(\X_i,\Z_i,\X_j,\Z_j),(\X_j, \Z_j) \big)|  > t\Big) 
    \leq  C_2 D \exp \Big( - \frac{ \alpha_2 nh^p t}{ 8 f_{\Z,max}}\int K^2 \Big),
    \label{eq:bound_major_uniform}
\end{align}
when $ t \leq 2 h^p  (\int K^2)^3 f^3_{\Z,max} /  C_K^4$.
It remains to bound $\PP(\sup_{\z\in\Zc} | n^{-1}\sum_{i=1}^n \bar g_i|>t/4) $.
Consider the family of functions $$\Fc:=\{ (\x_1,\z_1)\in \Rb\times \Zc\mapsto \frac{h^{p}}{4C_K^2}\EE[g_{\z}(\x_1,\z_1,\X,\Z)],\, \z\in \Zc \} .$$ 
This family of functions is bounded is one and its variance is less than
$\bar\sigma^2:=h^p\int K^2 f_{\z,max}^3 \big(  \int |K|\big)^2$.
Apply Propositions 9 and 10 in~\cite{fermanian2015single} that is coming from~\cite{einmahl2005uniform}: for some universal constants $A_1$ and $A_2$, some constant $A_{\bar g}$ that depends on $K$ and $f_{\z,max}$ (see Proposition 1 in~\cite{einmahl2005uniform}) and for every $x>0$,
\begin{eqnarray*}
\lefteqn{ \PP\Big(\sup_{\z\in\Zc} \frac{h^{p}}{4C_K^2} |\sum_{i=1}^n
 \EE[g_{\z}(\X_i,\Z_i,\X,\Z) | \X_i,\Z_i] | > A_1\big( x + A_{\bar g} n^{1/2}\bar\sigma \ln(1/\bar\sigma) \big)  \Big) }\\
&\leq & 2 \Big(\exp\big( -\frac{A_2 x^2}{n\bar\sigma^2}  \big) + \exp(-A_2 x) \Big),\;\text{or}
\hspace{3cm}
\end{eqnarray*}
$$\PP\Big(\sup_{\z\in\Zc} \frac{1}{n} |\sum_{i=1}^n
 \bar g_i | > 4A_1 C_K^2\big(  x - \frac{A_{\bar g}\bar\sigma}{n^{1/2}h^{p}} \ln(\bar\sigma) \big)  \Big) \leq 2 \exp\big( -\frac{A_2 n h^{2p}x^2}{\bar\sigma^2}  \big) 
 + 2\exp(-A_2 n h^p x).$$
For any positive $t$ s.t. $ 4 A_1 C_K^2 (n-1)A_{\bar g}\bar\sigma  \ln(1/\bar\sigma)  < n^{3/2}h^{p}t/8$, note that we can find a real $x> th^p / (16 C_K^2 A_1)$.
Then, we have
\begin{equation}
 \PP\Big(\sup_{\z\in\Zc} \frac{(n-1)}{n^2}|\sum_{i=1}^n \bar g_i | > \frac{t}{4}\Big) \leq 
2 \exp\big( -\frac{A_2 n h^p t^2 C_K^{-4}}{16^2 A_1^2  \int K^2 f_{\z,max}^3 (\int |K|)^2}  \big) + 2\exp(-\frac{A_2 n h^p t}{16 C_K^2 A_1} ) . 
\label{ineg_unif_barg_i}
\end{equation} 

\mds

Therefore, for such $t$, we obtain from~(\ref{ineg_unif_barg_i}) and~(\ref{eq:bound_major_uniform}) that
\begin{eqnarray*}
\lefteqn{ \PP \Big( \sup_{\z\in\Zc} | \sum_{1 \leq i \neq j \leq n}
    \big( S_{i,j}(\z) - \EE[S_{i,j}(\z)]  \big) |
    \geq t \Big)  \leq    C_2 D\exp \bigg( - \frac{ \alpha_2 nh^p t}{ 8(\int K^2) f_{\Z,max}} \bigg) }\\
    &+& 
    2 \exp\big( -\frac{A_2 n h^p t^2 C_K^{-4}}{15^2 A_1^2  \int K^2 f_{\z,max}^3 (\int |K|)^2}  \big) + 2\exp(-\frac{A_2 n h^p t}{15 C_K^2 A_1} ) . 
\end{eqnarray*}    
for sufficiently large integers $n$. The latter inequality,~(\ref{proba_delta_n_uniform}) and~(\ref{upper_bound_ES12_uniform}) yield the exponential upper bound. $\;\;\Box$

\medskip

\subsection{Proof of Proposition \ref{prop:consistency_hatTau}}
\label{proof:prop:consistency_hatTau}

Note that
$\tau_{1,2|\Z=\z} = \EE \big[ g_k(\X_{1} , \X_{2}) \big| \Z_1 = \z, \Z_2 = \z \big]$ for every $k=1,2,3$,
and that our estimators with the weights (\ref{def:weights_w_in}) can be written as
$\hat \tau_{1,2|\Z=\z}^{(k)} := U_{n} (g_k) \, / \, \{ U_{n} (1) + \epsilon_n\}$,
where
$$U_{n} (g) := \frac{1}{n(n-1) }
\sum_{1\leq i\neq j \leq n} g(\X_{i} , \X_{j})
\frac{K_{h} (\z - \Z_{i}) K_{h} (\z- \Z_{j})}{\EE[K_{h} (\z-\Z)]^{2}}=: \frac{1}{n(n-1)} \sum_{1\leq i\neq j \leq n} g_{i,j},$$
for any measurable bounded function $g$, with the residual diagonal term
$ \epsilon_n := \sum_{i=1}^n K_h^2(\z-\Z_i)/ \{n(n-1) \EE[K_{h} (\z-\Z)]^{2}\}$.
By Bochner's lemma  (see Bosq and Lecoutre~\cite{bosqlecoutre1987}),
$\epsilon_n$ is $O_P((nh^p)^{-1})$, and it will be negligible compared to $U_n(1)$.
Since the reasoning will be exactly the same for every estimator $\tau^{(k)}_{1,2|\z}$, i.e. for every function $g_k$, $k=1,2,3$, 
we omit the sub-index $k$. Then, the functions $g_k$ will be simply denoted by $g$. 

\mds

The expectation of our U-statistics is
\begin{eqnarray*}
\lefteqn{     \EE \big[ U_{n} (g) \big]
    :=  \EE \big[ g(\X_{1} , \X_{2})
    K_{h} (\z - \Z_{1}) K_{h} (\z - \Z_{2}) \big] /\EE[K_{h} (\z - \Z)]^{2} }\\
    &= &\int g(\x_{1} , \x_{2})
    K(\t_1) K(\t_2) f_{\X, \Z}(\x_1, \z + h \t_1) f_{\X,\Z}(\x_2, \z + h \t_2) d\x_1 \,d\x_2 \,d\t_1 \,d\t_2 /\EE[K_{h} (\z - \Z)]^{2} \\
    &\to & \frac{1}{f_\Z^2(\z)} \int g(\x_1, \x_2) f_{\X, \Z}(\x_1, \z) f_{\X, \Z}(\x_2, \z) d\x_1 d\x_2
    = \EE \big[ g(\X_{1} , \X_{2}) \big| \Z_1 = \z, \Z_2 = \z \big],
\end{eqnarray*}
applying Bochner's lemma to
$\z\mapsto \int g(\x_1,\x_2) f_{\X|\Z=\z}(\x_1)f_{\X|\Z=\z}(\x_2)\, d\x_1 \,d\x_2=\tau_{1,2|\Z=\z}$, that is a continuous function by assumption.

\mds

Set $\theta_n:=\EE[ U_n(g)]$, $g^*(\x_1, \x_2) := (g(\x_1, \x_2) + g(\x_2, \x_1))/2$ and $g^*_{i,j}=(g_{i,j}+g_{j,i})/2$ for every $(i,j)$, $i\neq j$. 
Note that $U_{n} (g) = U_{n} (g^*)$.
Since $g^*$ is symmetrical, the Hájek projection $\hat U_{n} (g^*)$ of $ U_{n} (g^*)$ satisfies
$\hat U_{n}(g^*)  := 2 \sum_{j=1}^n  \EE[g^*_{0,j} | \X_j,\Z_j]/n - \theta_n.$
Note that $\EE[\hat U_{n}(g^*)]=\theta_n= \tau_{1,2|\Z=\z}+o_P(1)$.
Since $ Var(\hat U_{n}(g^* ) = 4 Var(\EE[g^*_{0,j} | \X_j,\Z_j])/n=O((n h^p)^{-1})$, then $\hat U_{n}(g^*)=\theta_n+o_P(1)=\tau_{1,2|\Z=\z}+o_P(1)$.

\mds

Moreover, using the notation $\bar g_{i,j} := g^*_{i,j}- \EE[g^*_{i,j} | \X_j,\Z_j] -\EE[g^*_{i,j} | \X_i,\Z_i]+ \theta_n$ for $1 \leq i \neq j \leq n$,
we have $U_{n}(g^*) -\hat U_{n}(g^*) = \sum_{1\leq i\neq j\leq n}  \bar g_{i,j}/ n(n-1)$.
By usual U-statistics calculations, it can be easily checked that
\begin{eqnarray*}
Var\big( U_{n} (g^*) - \hat U_{n} (g^*)\big) &=&
\frac{1}{n^2(n-1)^2} \sum_{1\leq i_1\neq j_1\leq n} \sum_{1\leq i_2\neq j_2\leq n}
\EE[ \bar g_{i_1,j_1}\bar g_{i_2,j_2}] =O\big( \frac{1}{n^2h^{2p}} \big).
\end{eqnarray*}
Indeed, when all indices $(i_1,i_2,j_1,j_2)$ are different, or when there is a single identity among them,
$\EE[ \bar g_{i_1,j_1}\bar g_{i_2,j_2}]$ is zero. The
first nonzero terms arise when there are two identities among the indices, i.e. $i_1=i_2$ and $j_1=j_2$ (or $i_1=j_2$ and $j_1=i_2$).
In the latter case, we get an upper bound as $O((nh^p)^{-2})$ when $f_\Z$ is continuous at $\z$, by usual changes of variable techniques and Bochner's Lemma. Then,
$ U_{n}(g^*)=\hat U_n(g^*) + o_P(1)= \tau_{1,2|\Z=\z}+o_P(1)$.
Note that $U_{n} (1) + \epsilon_n$ tends to one in probability (Bochner's lemma).
As a consequence, $\hat \tau_{1,2|\Z=\z} = U_{n} (g^*) \, / \, (U_{n} (1)+\epsilon_n)$ tends to $\tau_{1,2|\Z=\z} / 1$ by the continuous mapping theorem.
$\Box$

\subsection{Proof of Proposition \ref{prop:unif_consistency_hatTau}}
\label{proof:prop:unif_consistency_hatTau}

Let us note that
\begin{align*}
    \tau_{1,2|\Z=\z} = \EE \big[ g_k(\X_{1} , \X_{2}) \big| \Z_1 = \z, \Z_2 = \z \big]
    = \int g_k(\x_{1} , \x_{2}) f_{\X|\Z=\z} (\x_1) f_{\X|\Z=\z} (\x_2) d\x_1 d\x_2
    = \frac{\phi_k(\z)}{f_{\Z}^2(\z)}, 
\end{align*}
where $\phi_k(\z) := \int g_k(\x_{1} , \x_{2})
f_{\X, \Z} (\x_1, \z) f_{\X, \Z} (\x_2, \z)
d\x_1 d\x_2$.
Also write
$\hat \tau_{1,2|\Z=\z}^{(k)} = \hat \phi_k(\z) / \hat f_{\Z}^2(\z)$, where $\hat \phi_k(\z) := n^{-2} 
\sum_{i,j=1}^n K_h(\Z_i - \z) K_h(\Z_j - \z) g_k(\X_i, \X_j)$
and $\hat f_{\Z}(\z) := n^{-1} \sum_{i=1}^n K_h(\Z_i - \z)$.
Therefore, we have
\begin{align*}
    \hat \tau_{1,2|\Z=\z}^{(k)} - \tau_{1,2|\Z=\z}
    = \frac{\hat \phi_k(\z) - \phi_k(\z)}{\hat f_{\Z}^2(\z)}
    - \tau_{1,2|\Z=\z} \frac{\hat f_{\Z}(\z) - f_{\Z}(\z)}{\hat f_{\Z}^2(\z)} \times
    \big( \hat f_{\Z}(\z) + f_{\Z}(\z) \big).
\end{align*}
By usual uniform consistency results (see for example Bosq and Lecoutre \cite{bosqlecoutre1987}),
$\sup_{\z \in \Zc} \big| \hat f_{\Z}(\z) - f_{\Z}(\z) \big| \to 0$ almost surely, as $n\to \infty$.
We deduce that $$\min_{\z\in \Zc}\hat f^2_{\Z}(\z)\geq f^2_{\Z,\min}/2,\, \text{ and} 
\max_{\z\in \Zc}| \hat f_{\Z}(\z)+ f_{\Z}(\z) | \leq 2\max_{\z\in\Zc} f_{\Z}(\z)\; \; \text{a.s.}$$
This means it is sufficient to prove the uniform strong consistency of $\hat \phi_k$ on $\Zc$, to obtain 
that $\sup_{\z \in \Zc} \big| \hat \tau_{1,2|\Z=\z}^{(k)} - \tau_{1,2|\Z=\z}^{(k)} \big|$ tends to zero a.s.

\mds

Note that, by Bochner's Lemma, $\sup_{\z \in \Zc} \big| \EE [\hat \phi_k(\z)] - \phi_k(\z) \big| \to 0$. 
Then, it remains to show that $\sup_{\z \in \Zc} \big| \hat \phi_k(\z) - \EE [\hat \phi_k(\z)] \big| \to 0$ almost surely.
Let $\rho_n > 0$ be such that we cover $\Zc$ by the union of $l_n$ open balls $B(\t_l, \rho_n)$, where $\t_1, \dots, \t_{l_n} \in \Rb^p$ and $l_n \in \NN^*$. Then
\begin{align*}
    \sup_{\z \in \Zc}
    \big| \hat \phi_k(\z) - \EE [\hat \phi_k(\z)] \big|
    \leq \sup_{l = 1, \dots l_n}
    \big| \hat \phi_k(\t_l) - \EE [\hat \phi_k(\t_l)] \big|
    + A_n,
\end{align*}
where $A_n := \sup_{l = 1, \dots l_n} \sup_{\z \in B(\t_l, \rho_n)}
\big| \hat \phi_k(\z) - \hat \phi_k(\t_l)
- (\EE [\hat \phi_k(\z)] - \EE [\hat \phi_k(\t_l)])\big|$.
For any index $l\in\{1,\ldots,l_n\}$ and any $\z\in B(\t_l, \rho_n)$, a first-order expansion yields
\begin{align*}  
    &\big| \hat \phi_k(\z) - \hat \phi_k(\t_l)
    - (\EE [\hat \phi_k(\z)] - \EE [\hat \phi_k(\t_l)]) \big| \\
    &= \bigg| \frac{1}{n(n-1)}
    \sum_{1 \leq i\neq j \leq n} g_k(\X_{i} , \X_{j})
    K_{h} (\z - \Z_{i}) K_{h} (\z- \Z_{j}) \\
    & \hspace{0.5cm} - \frac{1}{n(n-1)} \sum_{1 \leq i\neq j \leq n} g_k(\X_{i} , \X_{j})
    K_{h} (\t_l - \Z_{i}) K_{h} (\t_l- \Z_{j}) \\
    &\hspace{0.5cm} - \Big( \EE \big[ g_k(\X_{1} , \X_{2})
    K_{h} (\z - \Z_{1}) K_{h} (\z- \Z_{2}) \big]
    - \EE \big[ g_k(\X_{i} , \X_{j})
    K_{h} (\t_l - \Z_{i}) K_{h} (\t_l- \Z_{j}) \big] \Big) \bigg| \\
    &\leq  \frac{C_{Lip, K}}{h^{2p+1}}  |\z - \t_l|
    \Big( \EE \big[ | g_k(\X_{1} , \X_{2}) | \big] + \frac{1}{n(n-1)}
    \sum_{1 \leq i\neq j \leq n} |g_k(\X_{i} , \X_{j})| \Big) \\
    &= O(  \frac{\rho_n}{h^{2p +1}} ) = o(1),
\end{align*}
for some constant $ C_{Lip, K}$ and by choosing $\rho_n = o(h_n^{2p +1})$. Actually, we can cover $\Zc$ in such a way that $l_n = O(h_n^{-p(2p + 1)})$. 
This is always possible because $\Zc$ is a bounded set in $\Rb^p$. The previous upper bound is uniform w.r.t. 
$l$ and $\z \in B(\t_l, \rho_n)$, proving $A_n = o(1)$ everywhere.

\mds

Now, for every $l=\leq  l_n$, apply Equation (\ref{eq:bound_bernestein_improved}) for every $\z = \t_l$.
For any $t>0$, this yields
$$\PP \bigg( \frac{1}{n(n-1)}\Big| \sum_{i \neq j}
    g^{(l)} \big( (\X_i, \Z_i), (\X_j, \Z_j) \big)
    - \EE \Big[ g^{(l)} \big( (\X_1, \Z_1), (\X_2, \Z_2) \big) \Big] \Big| > t \bigg) 
    \leq \exp \Big( - \frac{C_0 nh_n^{2p} t^2 }{ C_1 + C_2 t} \Big),$$
for some positive constants $C_0,C_1,C_2$, by setting
$$g^{(l)} \big( (\X_i, \Z_i), (\X_j, \Z_j) \big)
:= g_k(\X_i , \X_j)
K_{h} (\t_l - \Z_i) K_{h} (\t_l - \Z_j).$$
Therefore, we deduce
$$    \PP \left( \sup_{l = 1, \dots l_n}
    \big| \hat \phi_k(\t_l) - \EE [\hat \phi_k(\t_l)] \big|  \geq t \right)
    \leq C_4 h_n^{-p(2p + 1)} \exp \Big( - \frac{C_0 nh_n^{2p} t^2 }{ C_1 + C_2 t} \Big),
$$
for some constant $C_4$.
Finally, applying Borel-Cantelli lemma,
$\sup_{\z \in \Zc} \big| \hat \phi_k(\z) - \EE [\hat \phi_k(\z)] \big|$ tends to zero a.s., proving the result. 
$\Box$

\medskip
    
\subsection{Proof of Proposition \ref{prop:asymptNorm_hatTau}}
\label{proof:prop:asymptNorm_hatTau}

By Markov's inequality, $\sum_{i=1}^n w_{i,n}^2(\z)=O_P((n h^p)^{-1})$ for any $\z$, that tends to zero. 
Then, by Slutsky's theorem, we get an asymptotic equivalence between the limiting laws of any $\hat\tau^{(k)}_{1,2|\z}$, $k=1,2,3$, and of their linearly transformed versions $\tilde\tau_{1,2|\z}$. Thus, we will prove the asymptotic normality of 
$\hat\tau^{(k)}_{1,2|\z}$ for some index $k=1,2,3$, simply denoted by $\hat \tau_{1,2|\z}$.

\mds
Let $g^*(\x_1, \x_2) := (g_k(\x_1, \x_2) + g_k(\x_2, \x_1))/2$ for some index $k=1,2,3$ (that will be implicit in the proof).
We now study the joint behavior of $(\hat \tau_{1,2|\Z=\z'_i} -
\tau_{1,2|\Z=\z'_i})_{i=1, \dots, n'}$.
We will extend Stute \cite{stute1991conditional}'s approach, in the case of multivariate conditioning variable $\z$ and studying the joint distribution of U-statistics at several conditioning points.
As in the proof of Proposition~\ref{prop:consistency_hatTau}, the estimator with the weights given by~(\ref{def:weights_w_in}) can be rewritten as
$\hat \tau_{1,2|\Z=\z'_i} := U_{n,i} (g^*) \, / \, (U_{n,i} (1)+ \epsilon_{n,i}),$
where
$$U_{n,i} (g) := \frac{1}{n(n-1) \EE[K_{h} (\z'_i - \Z)]^{2}}
\sum_{j_1, j_2 = 1,j_1\neq j_2}^{n} g(\X_{j_1} , \X_{j_2})
K_{h} (\z'_i - \Z_{j_1}) K_{h} (\z'_i - \Z_{j_2}),$$
for any bounded measurable function $g: \Rb^4 \to \Rb$.
Moreover, $\sup_{i=1,\ldots,n'}|\epsilon_{n,i}| = O_P(n^{-1}h^{-p})$.
By a limited expansion of $f_{\X,\Z}$ w.r.t. its second argument, and under
Assumption~\ref{assumpt:f_XZ_Holder}, we easily check that
$\EE \big[U_{n,i} (g) \big]  = \tau_{1,2|\Z = \z'_i} + r_{n,i}$,
where $|r_{n,i}| \leq C_0 h_n^\alpha/f^2_\Z(\z'_i)$, for some constant $C_0$ that is independent of $i$.

\mds

Now, we prove the joint asymptotic normality of $\big( U_{n,i} (g) \big)_{i=1, \dots, n'}$.
The Hájek projection $\hat U_{n,i} (g)$ of $U_{n,i} (g)$ satisfies
$\hat U_{n,i}(g) := 2 \sum_{j=1}^n  g_{n,i} \big(\X_{j}, \Z_j \big)/n - \theta_n$, where $\theta_n := \EE\big[U_{n,i} (g) \big]$ and
\begin{align*}
    &g_{n,i}(\x, \z) :=  K_{h} (\z'_i - \z) \EE \big[ g(\X , \x)
    K_{h} (\z'_i - \Z)  \big]
    \, / \, \EE[K_{h} (\z'_i - \Z)]^{2}.
\end{align*}

\begin{lemma}
    Under the assumptions of Proposition~\ref{prop:asymptNorm_hatTau}, for any measurable bounded function~$g$,
    \begin{align*}
        (n h^p)^{1/2} \Big( \hat U_{n,i}(g) - \EE \big[U_{n,i} (g) \big] \Big)_{i=1, \dots, n'}
        \indistrto \Nc(0, M_\infty (g)), \text{ as } n \to \infty,
    \end{align*}
    where, for $1 \leq i,j \leq n'$,
    \begin{align*}
        [M_\infty(g)]_{i, j} &:= \frac{ 4 \int K^2 \1_{ \{ \z'_i = \z'_j \} }  } { f_\Z(\z'_i)}
        \int g \big(\x_1, \x) g \big(\x_2, \x)
        f_{\X | \Z = \z'_i}(\x) f_{\X | \Z = \z'_i}(\x_1) f_{\X | \Z = \z'_i}(\x_2) d\x \, d\x_1 \, d\x_2.
    \end{align*}
    \label{lemma:Hajek_proj_U}
\end{lemma}
This lemma is proved in~\ref{proof:lemma:Hajek_proj_U}.
Similarly as in the proof of Lemma 2.2 in Stute \cite{stute1991conditional}, for every $i=1, \dots, n'$ and every bounded symmetrical measurable function $g$, we have
$(n h^p)^{1/2} Var \big[ \hat U_{n,i}(g) - U_{n,i}(g) \big] = o(1)$, which implies
\begin{align*}
    (n h^p)^{1/2} \Big( U_{n,i}(g) - \EE \big[U_{n,i} (g) \big] \Big)_{i=1, \dots, n'}
    \indistrto \Nc(0, M_\infty(g)), \text{ as } n \to \infty.
\end{align*}

Considering two measurable bounded functions $g_1$ and $g_2$, we have $U_{n,i}(c_1 g_1 + c_2 g_2)=c_1 U_{n,i}(g_1) + c_2 U_{n,i}(g_2)$ for every numbers $c_1, c_2$.
By the Cramér-Wold device, we check that
\begin{align*}
    &(n h^p)^{1/2} \bigg(
    \Big( U_{n,i}(g_1) - \EE \big[U_{n,i} (g_1) \big] \Big)_{i=1, \dots, n'},
    \Big( U_{n,i}(g_2) - \EE \big[U_{n,i} (g_2) \big] \Big)_{i=1, \dots, n'}
    \bigg) \\
    &\hspace{5cm} \indistrto \Nc \left(0,
    \begin{bmatrix}
        M_\infty(g_1)      & M_\infty(g_1, g_2) \\
        M_\infty(g_1, g_2) & M_\infty(g_2)
    \end{bmatrix}
    \right),
\end{align*}
as $n \to \infty$, where
$$    [M_\infty(g_1, g_2)]_{i, j} := \frac{ 4 \int K^2 \1_{ \{ \z'_i = \z'_j \} }
    }{ f_\Z(\z'_i)} \int g_1 \big(\x_1, \x) g_2 \big(\x_2, \x)  
    f_{\X | \Z = \z'_i}(\x) f_{\X | \Z = \z'_i}(\x_1) f_{\X | \Z = \z'_i}(\x_2) d\x \, d\x_1 \, d\x_2.
$$
Set $\tilde \tau_{1,2|\Z=\z'_i} := U_{n,i} (g^*) \, / \, U_{n,i} (1)$. Since $(nh_n^{p})^{1/2}\big(\hat\tau_{1,2|\Z=\z'_i} - \tilde \tau_{1,2|\Z=\z'_i}\big) =O_P\big((n h_n^{p})^{1/2}\epsilon_{n,i} \big)$ is $o_P(1),$
it is sufficient to establish the asymptotic law of
$(nh_n^{p})^{1/2}\big(\tilde \tau_{1,2|\Z=\z'_i} - \tau_{1,2|\Z=\z'_i}\big)$.
Since $\EE[U_{n,i} (1)] =1 + o((n h^p)^{-1/2})$ and
$\EE[U_{n,i} (g^*)] = \tau_{1,2|\Z=\z'_i} + o((nh_n^{p})^{-1/2})$, we get
\begin{align*}
    &(n h^p)^{1/2} \bigg(
    \Big( U_{n,i}(g^*) - \tau_{1,2|\Z=\z'_i} \Big)_{i=1, \dots, n'},
    \Big( U_{n,i}(1) - 1 \Big)_{i=1, \dots, n'}
    \bigg) \\
    &\hspace{5cm} \indistrto \Nc \left(0,
    \begin{bmatrix}
        M_\infty(g^*)    & M_\infty(g^*, 1) \\
        M_\infty(g^*, 1) & M_\infty(1)
    \end{bmatrix}
    \right), \text{ as } n \to \infty.
\end{align*}
Now apply the Delta-method with the function $\rho(\x, \y) := \x / \y$ where $\x$ and $\y$ are real-valued vectors of size $n'$ and the division has to be understood component-wise.
The Jacobian of $\rho$ is given by the $n' \times 2 n'$ matrix
$$J_\rho(\x, \y) =
\begin{bmatrix} Diag \big(y_1^{-1}, \dots y_{n'}^{-1} \big) \, , \,
Diag \big(-x_1 y_1^{-2}, \dots -x_{n'} y_{n'}^{-2} \big)
\end{bmatrix},$$
where, for any vector $v$ of size $n'$, $Diag(v)$ is the diagonal matrix whose diagonal elements are the $v_i$, with $i=1, \dots, n'$.
We deduce $(n h^p)^{1/2} \left( \tilde  \tau_{1,2|\Z=\z'_i} - \tau_{1,2|\Z=\z'_i} \right)_{i=1, \dots, n'}
    \indistrto \Nc (0,  \HH), \text{ as } n \to \infty$, setting
\begin{align*}
    \HH := J_\rho(\vec{\tau}, \e)
    \begin{bmatrix}
        M_\infty(g^*)    & M_\infty(g^*, 1) \\
        M_\infty(g^*, 1) & M_\infty(1)
    \end{bmatrix}
    J_\rho(\vec{\tau}, \e)^T,
\end{align*}
where $\vec{\tau} =  \left( \tau_{1,2|\Z=\z'_i} \right)_{i=1, \dots, n'}$ and $\e$ is the vector of size $n'$ whose all components are equal to $1$.
Thus, we have $J_\rho(\vec{\tau}, \e) = \begin{bmatrix} Id_{n'} , -Diag(\vec{\tau}) \end{bmatrix}$, denoting by $Id_{n'}$ the identity matrix of size $n'$
and by $Diag(\vec{\tau})$ the diagonal matrix of size $n'$ whose diagonal elements are the $\tau_{1,2|\z'_i}$, for $i=1, \dots, n'$.
To be specific, we get
\begin{align*}
    \HH = M_\infty(g^*) - Diag(\vec{\tau}) M_\infty(g^*, 1)
    - M_\infty(g^*, 1) Diag(\vec{\tau})
    + Diag(\vec{\tau}) M_\infty(1) Diag(\vec{\tau}).
\end{align*}
For $i,j$ in $\{1,\ldots,n'\}$ and using the symmetry of the function $g^*$, we obtain
\begin{align*}
    [M_\infty(g^*)]_{i, j}
    &=  4 \int K^2 \1_{ \{ \z'_i = \z'_j \} }
    \EE[ g^*(\X_1,\X)g^*(\X_2,\X)| \Z = \Z_1 = \Z_2 = \z'_i]/f_\Z(\z'_i),
\end{align*}
\begin{eqnarray*}
\lefteqn{    [Diag(\vec{\tau}) M_\infty(g^*, 1)]_{i, j}
    = 4\tau_{1,2|\Z=\z'_i} \int K^2 \1_{ \{ \z'_i = \z'_j \} }
    \EE[ g^*(\X_1,\X)| \Z = \Z_1 = \z'_i]/f_\Z(\z'_i) }\\
    &=&  4 \int K^2 \1_{ \{ \z'_i = \z'_j \} } \tau_{1,2|\Z=\z'_i}^2/ f_\Z(\z'_i)
    = [M_\infty(g^*, 1) Diag(\vec{\tau})]_{i,j}
    = [Diag(\vec{\tau}) M_\infty(1) Diag(\vec{\tau})]_{i, j}.
\end{eqnarray*}
As a consequence, we obtain
\begin{align*}
    [\HH]_{i, j} = \frac{ 4 \int K^2 \1_{ \{ \z'_i = \z'_j \} }
    }{ f_\Z(\z'_{i})} \Big(
    \EE[g^*(\X_1,\X)g^*(\X_2,\X)  | \Z = \Z_1 = \Z_2 = \z'_i]- \tau_{1,2|\Z=\z'_i}^2
    \Big).\;\; \Box
\end{align*}

\medskip

\subsection{Proof of Lemma \ref{lemma:Hajek_proj_U} }
\label{proof:lemma:Hajek_proj_U}

Let us first evaluate the variance-covariance matrix
$M_{n,n'} := [Cov(\hat U_{n,i}, \hat U_{n,j})]_{1 \leq i, j \leq n'}$.
Note that $\EE \big[ g_{n,i} (\X_{j}, \Z_j ) \big]
= \EE \big[\hat U_{n,i} \big]
= \EE \big[ U_{n,i} (g) \big],$
and that
\begin{align*}
    &\Big( (n h^p)^{1/2} \big( \hat U_{n,i} - \EE [U_{n,i} (g) ] \big) \Big)_{i=1, \dots, n'} =
    \frac{2 h^{p/2}}{n^{1/2}} \sum_{j=1}^n
    \big( g_{n,i} (\X_{j}, \Z_j ) - \EE [U_{n,i} (g) ] \big)_{i=1, \dots, n'},
\end{align*}
that is a sum of independent vectors.
Thus,
$Cov(\hat U_{n,i}, \hat U_{n,j}) = 4 n^{-1} Cov\big( g_{n,i} \big(\X, \Z \big) , g_{n,j} \big(\X, \Z \big)  \Big) ,$ for every $i,j$ in $\{1,\ldots,n'\}$, and
\begin{eqnarray*}
\lefteqn{
\EE \big[ g_{n,i} (\X, \Z ) g_{n,j} (\X, \Z ) \big]   }\\
&=&
\int K_{h} (\z'_i - \z) K_{h} (\z'_j - \z) \frac{\EE \big[ g(\X , \x)
     K_{h} (\z'_i - \Z) \big]\EE \big[ g(\X , \x)
     K_{h} (\z'_j - \Z) \big]}
    {\EE[K_{h} (\z'_i - \Z)]^{2}\EE[K_{h} (\z'_j - \Z)]^{2}}
    f_{\X,\Z}(\x, \z) d\x \, d\z \\
&\sim & \frac{1}{h^{p} f_\Z^2(\z'_i)f_\Z^2(\z'_j) }
    \int g \big(\x_1, \x) g \big(\x_2, \x) K_{h} (\z'_i - \z) K_{h} (\z'_j - \z)
    K_{h} (\z'_i - \w_1) K_{h} (\z'_j - \w_2) \\
    & \times & f_{\X,\Z}(\x, \z) f_{\X,\Z}(\x_1, \w_1) f_{\X,\Z}(\x_2, \w_2)
    d\x \, d\z \, d\x_1 \, d\w_1 \, d\x_2 \, d\w_2 \\
    &\sim & \frac{1}{h^{p} f_\Z^2(\z'_i)f_\Z^2(\z'_j) }
    \int g \big(\x_1, \x) g \big(\x_2, \x)
    K(\u_1) K(\u_2) K (\u) K (\frac{\z'_j-\z'_i}{h} + \u) f_{\X,\Z}(\x, \z'_i - h \u) \\
    &  \times & f_{\X,\Z}(\x_1, \z'_i - h \u_1) f_{\X,\Z}(\x_2, \z'_j - h \u_2)
     d\x \, d\u \, d\x_1 \, d\u_1 \, d\x_2 \, d\u_2.
\end{eqnarray*}
If $i\neq j$ and $K$ is compactly supported, the latter term is zero when $n$ is sufficiently large, and
$ Cov(\hat U_{n,i}, \hat U_{n,j}) = - 4 n^{-1} \EE[U_{n,i}] \EE[U_{n,j}]
\sim - 4 n^{-1} \tau_{1,2|\Z=\z'_i} \tau_{1,2|\Z=\z'_j}=o\big((nh^p)^{-1}\big).$
Otherwise, $i=j$ and, as $\EE \Big[ g_{n,i} \big(\X_{1}, \Z_1 \big) \Big]  = O(1)$, we have
\begin{align*}
    Var & \Big( \big( g_{n,i} (\X, \Z ) \big)^2 \Big)
    \sim \frac{1}{h^p f_\Z^4(\z'_i)}
    \int g \big(\x_1, \x) g \big(\x_2, \x)
    K(\u_1) K(\u_2) K^2 (\u) f_{\X,\Z}(\x, \z'_i-h\u) \\
    & \hspace{2cm} \times f_{\X,\Z}(\x_1, \z'_i - h\u_1) f_{\X,\Z}(\x_2, \z'_i-h\u_2)
    \, d\x \, d\u \, d\x_1 \, d\u_1 \, d\x_2 \, d\u_2 \\
    &\sim \frac{\int K^2}{h^p f_\Z(\z'_i)}
    \int g \big(\x_1, \x) g \big(\x_2, \x)
    f_{\X | \Z = \z'_i}(\x) f_{\X | \Z = \z'_{i}}(\x_1) f_{\X | \Z = \z'_i}(\x_2) \,d\x \, d\x_1 \, d\x_2,
\end{align*}
by Bochner's lemma.
We have proved that, for every $i,j \in \{1,\ldots,n'\}$,
$$ n h^p [M_{n,n'}]_{i, j} \to
    \frac{ 4 \int K^2 \1_{ \{ \z'_i = \z'_j \} }} { f_\Z(\z'_i)}
    \int g \big(\x_1, \x) g \big(\x_2, \x)
    f_{\X | \Z = \z'_i}(\x) f_{\X | \Z = \z'_i}(\x_1) f_{\X | \Z = \z'_i}(\x_2) \,d\x \, d\x_1 \, d\x_2,$$
as $n \to \infty$. Therefore, $n h^p M_{n,n'}$ tends to $M_\infty$.

\mds

We now verify Lyapunov's condition with third-order moments, so that the usual multivariate central limit theorem would apply.
It is then sufficient to show that
\begin{equation}
    \Big(\frac{ h^{p/2}}{n^{1/2}} \Big)^3  \sum_{j=1}^n \EE \Big[ \big| g_{n,i} (\X_{j}, \Z_j ) - \,
    \EE [U_{n,i} (g) ] \big|^3 \Big]
    = o(1).
    \label{cond:CLT_Lyapunov}
\end{equation}

For any $j=1,\ldots,n$, we have
\begin{align*}
    &\EE \Big[ \big| g_{n,i} (\X_{j}, \Z_j ) - \, \EE [U_{n,i} (g) ] \big|^3 \Big] \\
    &\sim \int \Big| \frac{1}{f_\Z^2(\z'_i)} \int g( \x_1, \x) K_{h} (\z'_i -  \z_1) K_{h} (\z'_i - \z)
    f_{\X, \Z}( \x_1,  \z_1) d \x_1\, d \z_1
    - \, \EE \big[U_{n,i} (g) \big]  \Big|^3 f_{\X, \Z}(\x, \z) d\x\, d\z.
\end{align*}
By the change of variable $ \z_1 = \z'_i - h \t_1$ and $\z = \z'_i - h \t$, we get
\begin{eqnarray*}
    \lefteqn{
    \EE \Big[ \big| g_{n,i} (\X_{j}, \Z_j ) - \, \EE [U_{n,i} (g) ] \big|^3 \Big]
    \sim  h^{-2p}\int \Big| \frac{1}{f_\Z^2(\z'_i)} \int g( \x_1, \x) K(\t_1) K(\t) f_{\X, \Z}(\x_1, \z'_i - h \t_1) d \x_1 \,d\t_1 }\\
    &-& h^p  \EE \big[U_{n,i} (g) \big]  \Big|^3 f_{\X, \Z}(\x, \z'_i - h \t) d\x\, d\t = O(h^{-2p}),\hspace{5cm}
\end{eqnarray*}
because of Bochner's lemma, under our assumptions.
Therefore, we have obtained
\begin{align*}
    \Big(\frac{h^{p/2}}{n^{1/2}} \Big)^3  \sum_{j=1}^n \EE \Big[ \big| g_{n,i} (\X_{j}, \Z_j ) - \,
    \EE [U_{n,i} (g) ] \big|^3 \Big]
    = O(h^{3p/2} n^{-3/2} n h^{-2p})
    = O((n h^{p})^{-1/2}) = o(1).
\end{align*}
Therefore, we have checked Lyapunov's condition and the result follows.
$\;\; \Box$

\end{document}